\documentclass[10pt]{amsart}
\usepackage{amsfonts,amsmath,mathtools}
\usepackage{multirow}
\usepackage{amssymb}
\usepackage{latexsym}
\usepackage{epsfig}
\usepackage{graphicx,color,graphics}
\usepackage{subcaption}

\newtheorem{lemma}{Lemma}[section]

\newtheorem{theorem}[lemma]{Theorem}
\newtheorem{proposition}[lemma]{Proposition}

\newtheorem{remark}[lemma]{Remark}

\newcommand{\N}{\ifmmode{{\Bbb N}}\else{\mbox{${\Bbb N}$}}\fi}
\newcommand{\R}{\ifmmode{{\Bbb R}}\else{\mbox{${\Bbb R}$}}\fi}

\renewcommand{\theequation}{\thesection.\arabic{equation}}
\newcommand{\halfscript}[2]{%
  \mathord{\hbox{
    \valign{%
      \vfil##\vfil\cr
      \hbox{$#1$}\cr
      \hbox{$\scriptstyle#2$}\cr
    }%
    \kern\scriptspace
  }}%
}

\begin{document}
\thispagestyle{empty}

\title[Rao--Nakra sandwich beam with a fractional dissipation]{Stability of the Rao--Nakra sandwich beam with a dissipation of
fractional derivative type: theoretical and numerical study}

\author{K. Ammari}
\address{LR Analysis and Control of PDEs, LR 22ES03, Department of Mathematics, Faculty of Sciences of Monastir, University of Monastir, 5019 Monastir, Tunisia}
\email{kais.ammari@fsm.rnu.tn}

\author{V. Komornik}
\thanks{The author was supported by the following grants: NSFC No. 11871348, CAPES: No. 88881.520205/2020-01,  
MATH AMSUD: 21-MATH-03.}
\address{D\'epartement de math\'ematique, Universit\'e de Strasbourg, 7 rue René Descartes, 67084 Strasbourg Cedex, France}
\email{komornik@math.unistra.fr}

\author{M. Sep\'{u}lveda}
\thanks{The author was supported by the following grants: Fondecyt-ANID project 1220869,
and ANID-Chile through Centro de Modelamiento  Matem\'atico  (FB210005)}
\address{DIM and CI2MA, Universidad del Concepci\`on, Concepci\`on, Chile}
\email{maursepu@udec.cl}

\author{O. Vera} 
\thanks{The author is partially financed by project Fondecyt 1191137 and UTA MAYOR 2022-2023, 4764-22.}
\address{Departamento de Matem\'{a}tica, Universidad
de Tarapac\'{a}, Av. 18 de Septiembre 2222, Arica, 
Chile}
\email{opverav@academicos.uta.cl}

\begin{abstract}
 This paper is devoted to the solution and stability  of a one-dimensional model depicting Rao--Nakra sandwich beams, incorporating damping terms characterized by fractional derivative types within the domain, specifically a generalized Caputo derivative with exponential weight. To address existence, uniqueness, stability, and numerical results, fractional derivatives are substituted by diffusion equations relative to a new independent variable, $\xi$, resulting in an augmented model with a dissipative semigroup operator. Polynomial decay of energy is achieved, with a decay rate depending on the fractional derivative parameters. Both the polynomial decay and its dependency on the parameters of the generalized Caputo derivative are numerically validated. To this end, an energy-conserving finite difference numerical scheme is employed.
\end{abstract}


 
 \maketitle
\tableofcontents

\section{Introduction}
\setcounter{equation}{0}

In 1999, Liu, Trogdon and Yong \cite{lty} developed the following general three-layer laminated beam--plate model: 
\begin{eqnarray}
\left\lbrace
\begin{array}{l}
\label{101}
\varrho_{1}h_{1}u_{\bf tt} - E_{1}h_{1}u_{xx} - \tau = 0,   \\
\varrho_{3}h_{3}v_{\bf tt} - E_{3}h_{3}u_{xx} + \tau = 0,   \\
\varrho h w_{\bf tt} + EIw_{xxxx} - G_{1}h_{1}(w_{x} + \phi_{1})_{x} - G_{3}h_{3}(w_{x} + \phi_{3})_{x} - h_{2}\tau_{x} = 0,   \\
\varrho_{1}I_{1}\phi_{1{\bf tt}} - E_{1}I_{1}\phi_{1xx} - \frac{h_{1}}{2}\tau + G_{1}h_{1}(w_{x} + \phi_{1}) = 0, \\
\varrho_{3}I_{3}\phi_{3{\bf tt}} - E_{3}I_{3}\phi_{3xx} - \frac{h_{3}}{2}\tau + G_{3}h_{3}(w_{x} + \phi_{3}) = 0.
\end{array}
\right.
\end{eqnarray}
The physical parameters $h_{i},$ $\varrho_{i},$ $E_{i},$ $G_{i},$ $I_{i}$ $>0$ are the thickness, density, Young's modulus, shear modulus and moments of inertia of the $i$-th layer for $i = 1,\,2,\,3,$ from bottom to top, respectively. In addition, $\varrho h = \varrho_{1}h_{1} + \varrho_{3}h_{3}$ and $EI = E_{1}I_{1} + E_{3}I_{3}.$ 

\medskip

The Rao--Nakra sandwich beam model \cite{rao} is the following:
\begin{eqnarray}
\label{102} \left\lbrace
\begin{array}{l}
\rho_{1}u_{\bf tt}(x,\,{\bf t}) - \vartheta u_{xx}(x,\,{\bf t}) - k(-u(x,\,{\bf t}) + v(x,\,{\bf t}) + \gamma w_{x}(x,\,{\bf t})) = 0,   \\
\rho_{2}v_{\bf tt}(x,\,{\bf t}) - \chi v_{xx}(x,\,{\bf t}) + k(-u(x,\,{\bf t}) + v(x,\,{\bf t}) + \gamma w_{x}(x,\,{\bf t})) = 0,   \\
\rho_{3}w_{\bf tt}(x,\,{\bf t}) + \zeta w_{xxxx}(x,\,{\bf t}) - k (-u(x,\,{\bf t}) + v(x,\,{\bf t}) + \gamma w_{x}(x,\,{\bf t}))_{x} = 0, \\
u(x,\,0) = u_{0}(x),\ u_{\bf t}(x,\,0) = u_{1}(x),   \\
v(x,\,0) = v_{0}(x),\ v_{\bf t}(x,\,0) = v_{1}(x),   \\
w(x,\,0) = w_{0}(x),\ w_{\bf t}(x,\,0) = w_{1}(x),
\end{array}
\right.
\end{eqnarray}
with $0<x<\ell$ and ${\bf t}>0.$ Here $u = u(x,\,{\bf t})$ and $v = v(x,\,{\bf t})$ are the longitudinal displacement. $w = w(x,\,{\bf t})$ is the transverse displacement of the beam. This system 
is derived from the system \eqref{101} when we consider the core material to be linearly elastic, that is, $\tau = 2G_{2}\gamma$ with the shear strain 
\begin{align*}
\gamma = \frac{1}{2h_{2}}(-u + v + \alpha w_{x})\quad {\rm and}\quad \alpha = h_{2} = \frac{1}{2}(h_{1} + h_{3}),
\end{align*}
where $k = \frac{G_{2}}{h_{2}},$ the shear modulus $G_{2} = \frac{E_{2}}{2(1 + \nu)},$ and $-1<\nu<\frac{1}{2}$ is the Poisson ratio.

\medskip

Many authors have been working on the system \eqref{102} from different points of view: see \cite{lrz, o1, o2, ra} and the references therein. 
In the past three decades there has been a growing interest by many reserchers for the study of fractional calculus  in different fields of sciences \cite{magin, tara, valerio}. Several aspects of engineering, applied sciences, and mathematical physics  benefitted  from this ascending wave of applications. Space sciences, fluid mechanics, porous media flows, viscoelastic and biological processes are but a few areas in which fractional order differential equations have become a favored tool to tread new path. To give some examples, fractional derivatives have been used to model frequency dependent damping behavior of many viscoelastic materials as well modeling many chemical processes. Many problems in several scientific applied areas, including the analysis of viscoelastic materials, heat conduction in materials with memory, electrodynamics with memory, signal processing, among others, can be modeled by fractional
differential calculus. Indeed, many investigations have
shown that models involving fractional derivatives are more
realistic to represent some natural phenomena than models involving classical derivatives. For more information we refer to \cite{KST, 15, SKM} and the references therein. For example, the wave equation with boundary fractional damping has been treated in \cite{15, mm} where  the strong stability and the lack of uniform stabilization were proved. Applications of this kind are also found in \cite{maryati, mu, po, ramos, vera} and the references therein. 

\medskip

In \cite{vera} the authors studied the Rao--Nakra system with a boundary dissipation of fractional derivative type. They established the polynomial stability of the system. 
In this article we are interested in studying the  polynomial stabilization to the system \eqref{103} 
from a theoretical and numerical point of view. In order to achieve this, we eliminate the dissipation given on the border and place the dissipation in the domain. The main idea to reach this result is to use an idea given in \cite{kais}. It is important to note that the imposed dissipation  is of a fractional derivative type. Our method is based on the fact that the input--output relationship in a certain  diffusion equation 
introduced by Mbodje in \cite{15} is realized by a
fractional derivative operator. 

\medskip

Numerical approximations of derivatives and fractional integrals have been extensively studied, developed, and refined. Comprehensive reviews on this topic are available in \cite{Cai-Li, li}. Many of these approaches involve direct approximations of the Riemann--Liouville, Caputo \cite{das,jesus,C1,C2,C3}, or Grünwald--Letnikov integral \cite{li}, using quadrature methods, finite differences, and truncated summations \cite{Ortigueira}. While some progress has been made in enhancing the convergence rate, such as through trapezoidal approximation \cite{li} or spline interpolation \cite{Ciesielski}, these methods do not inherently preserve the energy conservation property of the Rao--Nakra model under study, nor do they maintain stability in the presence of dissipative terms. In this work, we address this challenge by employing a suitable combination of the $\beta$-Newmark  method  \cite{Newmark}, coupled with a Crank--Nicolson method for the Mbodje diffusion equation, to achieve our objective.

Without lost of generality, we will consider the following system for $\rho_{1} = \rho_{2} = \rho_{3}=1$:
\begin{eqnarray}
\label{103} \left\lbrace
\begin{array}{l}
u_{\bf tt} - \vartheta u_{xx} - k(-u + v + \gamma w_{x}) + \partial_{\bf t}^{\alpha,\,\eta}u = 0,   \, (0,\ell) \times (0,+\infty),\\
v_{\bf tt} - \chi v_{xx} + k(-u + v + \gamma w_{x}) + \partial_{\bf t}^{\alpha,\,\eta}v = 0,   \, (0,\ell) \times (0,+\infty),\\
w_{\bf tt} + \zeta w_{xxxx} - k (-u + v + \gamma w_{x})_{x} + \partial_{\bf t}^{\alpha,\,\eta}w = 0,  \, (0,\ell) \times (0,+\infty),\\
u(x,\,0)=u_{0}(x),\;u_{\bf t}(x,\,0) = u_{1}(x),   \, (0,\ell),\\
v(x,\,0)=v_{0}(x),\;v_{\bf t}(x,\,0) = v_{1}(x),   \, (0,\ell),\\
w(x,\,0) = w_{0}(x),\;w_{\bf t}(x,\,0) = w_{1}(x), \, (0,\ell),\\
u(0,\,{\bf t}) = u(\ell,\,{\bf t}) = 0,\ v(0,\,{\bf t}) = v(\ell,\,{\bf t}) = 0, \, (0,+\infty),\\ 
w(0,\,{\bf t}) = w(\ell,\,{\bf t}) = w_{x}(0,\,{\bf t}) = w_{x}(\ell,\,{\bf t}) = 0, \, (0,+\infty),
\end{array}
\right.
\end{eqnarray}
with  $0 < x < \ell,$\ ${\bf t}>0$, and real-valued functions $u = u(x,\,{\bf t})$ and $v = v(x,\,{\bf t}),$ $w = w(x,\,{\bf t})$.

\medskip

The rest of the paper is divided into five sections. In Section \ref{2} we show that the  system \eqref{103} may be replaced by an augmented system \eqref{213} obtained by coupling an equation with a suitable diffusion, and we study the energy functional associated to system. In Section \ref{3}, we establish the existence and uniqueness of solutions of the system \eqref{213}; for this we use \cite{kaisbis,kais}. In section \ref{4} we prove the strong stability, the lack of uniform stabilization and the polynomial stability of the system \eqref{213}. In Section \ref{5}, we give a numerical study of the polynomial stability. 

\medskip

Throughout this paper, $C$ is a generic constant, not necessarily the same  at each occasion (it may change from  line to line) and depending on the indicated quantities.

\renewcommand{\theequation}{\thesection.\arabic{equation}}
\setcounter{equation}{0} \label{2}
\section{Augmented Model}
This section is devoted to the reformulation of the
model \eqref{103} with boundary conditions \eqref{103}$_{7,\,8}$ as an
augmented system. 
First we provide a brief review of fractional calculus. There are many definitions of fractional derivatives \cite{das}, among them those of Riemann--Liouville  and Caputo  are the most widely used \cite{jesus}. The latter has the same Laplace transform as the integer order derivative, so it is widely used in control theory. In this paper the fractional derivative damping force is regarded as a control force to study the properties of free damped vibrations of the system, hence the Caputo definition \cite{C1,C2,C3} is used here. 
Let $0 < \alpha < 1$.  The Caputo fractional integral of order
$\alpha$ is defined by the formula
\begin{equation}
\label{201}I^{\alpha}f({\bf t}) =
\frac{1}{\Gamma(\alpha)}\int_{0}^{\bf t}({\bf t} - s)^{\alpha - 1}
f(s)ds,
\end{equation}
where $\Gamma$ is the  Gamma function, $f \in
L^{1}(0,+\infty)$ and  ${\bf t} > 0.$ \\
\\
The Caputo fractional derivative operator of order $\alpha$ is
defined by the expression
\begin{equation}
\label{202} D^{\alpha}f({\bf t}) = I^{1 - \alpha}D f({\bf t}) :=
\frac{1}{\Gamma(1 - \alpha)}\int_{0}^{\bf t}({\bf t} - s)^{-\alpha}
f^\prime (s) ds.
\end{equation}
 
The Caputo definition of fractional derivative possesses a very simple but interesting  interpretation: if the function $f({\bf t})$
represents the strain history within a viscoelastic material whose
relaxation function is $[\Gamma(1 - \alpha){\bf t}^{\alpha}]^{-1}$, then
the material will experience at any time ${\bf t}$ a total stress given by
the expression $D^{\alpha}f({\bf t}).$ 
Also, it easy to show that  $D^{\alpha}$ is a left inverse of
$I^{\alpha},$  but in general it is not a right inverse. More
precisely, we have
\begin{eqnarray*}
D^{\alpha}I^{\alpha}f =  f, \qquad I^{\alpha}D^{\alpha}f({\bf t}) =  f({\bf t})
- f(0).
\end{eqnarray*}
See \cite{SKM} for the proof of above equalities and for more properties of fractional calculus. 

\medskip

In this work, we consider a slightly different version of
\eqref{201} and \eqref{204}.  In  \cite{CM}, Choi and
MacCamy introduced the following definition of fractional
integro-differential operators  with an exponential weight. Given  $0 <
\alpha < 1$ and $ \eta \ge 0$,  the exponential  fractional integral of
order $\alpha$ is defined by
\begin{equation}
\label{203}I^{\alpha,\ \eta}f({\bf t}) =
\frac{1}{\Gamma(\alpha)}\int_{0}^{\bf t}e^{-\eta({\bf t} - s)} ({\bf t} -
s)^{\alpha - 1}f(s)ds,
\end{equation}
and the exponential fractional derivative operator of order $\alpha$ is
defined by
\begin{equation}
\label{204}\partial_{\bf t}^{\alpha,\ \eta}f({\bf t}) = \frac{1}{\Gamma(1 - \alpha)}\int_{0}^{\bf t}e^{-\eta\,({\bf t} - s)}({\bf t} -
s)^{-\alpha} f^\prime (s) ds.
\end{equation}
Note that $\partial_{\bf t}^{\alpha,\ \eta}f({\bf t}) = I^{1 - \alpha,\,\eta}f'({\bf t}).$ 

\medskip

We will need later the following results:
 
\begin{theorem}\cite{15}
\label{theorem51} Let $\mu$ be the function
\begin{eqnarray}
\label{205}\mu(\xi) = |\xi|^{(2\alpha - 1)/2},\qquad
\xi\in\mathbb{R},\qquad 0<\alpha<1.
\end{eqnarray}
Then the relation between the {\it Input} $U$ and the {\it Output}
$O$ is given by the following system:
\begin{eqnarray}
\label{206}&  & \varphi_{\bf t}({\bf t},\,\xi) + |\xi|^{2}\varphi({\bf t},\,\xi) =
\mu(\xi)U({\bf t}),\qquad \xi\in\mathbb{R},\qquad {\bf t}>0, \\
\label{207}&  & \varphi(0,\,\xi) = 0,\\
\label{208}&  & O =
\pi^{-1}\sin(\alpha\pi)\int_{\mathbb{R}}\mu(\xi)\varphi({\bf t},\,\xi)
d\xi.
\end{eqnarray}
This implies that
\begin{eqnarray}
\label{209}O = I^{1 - \alpha}U,
\end{eqnarray}
where $U\in C([0,\,+\infty)).$
\end{theorem}

\begin{lemma}\cite[Lemma 2.1]{achouri}
\label{lemma22} If $\lambda \in D = \{\lambda\in \mathbb{C}:\ Re\,\lambda + \eta >0 \text{ or } Im\,\lambda \neq 0\}$, 
then
\begin{eqnarray*}
\int_{\mathbb{R}}\frac{\mu^{2}(\xi)}{\xi^{2} + \eta + \lambda} \, d\xi
= \frac{\pi}{\sin(\alpha\pi)}(\eta + \lambda)^{\alpha - 1}.
\end{eqnarray*}
\end{lemma}

Our strategy requires the elimination of the fractional derivatives in time from the domain condition in the system \eqref{103}. To this, setting $\mu(\xi) = |\xi|^{(2\alpha - 1)/2},$ $\xi\in\mathbb{R},$ $\mathfrak{C} = \pi^{-1}\sin(\alpha\,\pi),$ and exploiting the technique from \cite{huang}, we  reformulate the system \eqref{103} by using Theorem
\ref{theorem51}.
This leads to  the following augmented
model, where $(x,\,t,\,\xi) \in (0,\,\ell) \times(0,\,+\infty) \times \mathbb{R}$ and we introduce the constant  $\mathfrak{C} := \pi^{-1}\sin(\alpha\,\pi)$:
\begin{eqnarray}
\label{213} \left\lbrace
\begin{array}{l}
u_{\bf tt}(x,\,{\bf t}) - \vartheta u_{xx}(x,\,{\bf t}) - k(-u(x,\,{\bf t}) + v(x,\,{\bf t}) + \gamma w_{x}(x,\,{\bf t})) \\
+\ \mathfrak{C}\int_{\mathbb{R}}\mu(\xi)\varphi(x,\,{\bf t},\,\xi)
d\xi = 0,   \\
\varphi_{\bf t}(x,\,{\bf t},\,\xi) + (\xi^{2} + \eta)\varphi(x,\,{\bf t},\,\xi) = \mu(\xi)u_{\bf t}(x,\,{\bf t}), \\
\\
v_{\bf tt}(x,\,{\bf t}) - \chi v_{xx}(x,\,{\bf t}) + k(-u(x,\,{\bf t}) + v(x,\,{\bf t}) + \gamma w_{x}(x,\,{\bf t})) \\
+\ \mathfrak{C}\int_{\mathbb{R}}\mu(\xi)\phi(x,\,{\bf t},\,\xi)
d\xi = 0,   \\
\phi_{\bf t}(x,\,{\bf t},\,\xi) + (\xi^{2} + \eta)\phi(x,\,{\bf t},\,\xi) = \mu(\xi)v_{\bf t}(x,\,{\bf t}), \\
\\
w_{\bf tt}(x,\,{\bf t}) + \zeta w_{xxxx}(x,\,{\bf t}) - k(-u(x,\,{\bf t}) + v(x,\,{\bf t}) + \gamma w_{x}(x,\,{\bf t}))_{x} \\
+\ \mathfrak{C}\int_{\mathbb{R}}\mu(\xi)\psi(x,\,{\bf t},\,\xi)
d\xi = 0, \\
\psi_{\bf t}(x,\,{\bf t},\,\xi) + (\xi^{2} + \eta)\psi(x,\,{\bf t},\,\xi) = \mu(\xi)w_{\bf t}(x,\,{\bf t}), \\
\\
u(x,\,0)=u_{0}(x),\;u_{\bf t}(x,\,0) = u_{1}(x),   \\
v(x,\,0)=v_{0}(x),\;v_{\bf t}(x,\,0) = v_{1}(x),   \\
w(x,\,0) = w_{0}(x),\;w_{\bf t}(x,\,0) = w_{1}(x),  \\
u(0,\,t) = u(\ell,\,{\bf t}) = 0,\ v(0,\,{\bf t}) = v(\ell,\,{\bf t}) = 0, \\ 
w(0,\,t) = w(\ell,\,{\bf t}) = w_{x}(0,\,{\bf t}) = w_{x}(\ell,\,{\bf t}) = 0  \\
\varphi(x,0,\xi) = 0,\ \phi(x,0,\xi) = 0,\ \psi(x,0,\,\xi) = 0.
\end{array}
\right.
\end{eqnarray}

\section{Setting of the Semigroup} \label{3}
\setcounter{equation}{0}
In this section, we discuss the existence and uniqueness of solutions for the coupled system \eqref{213} by the  semigroup theory \cite{Pazy, prus}.

\medskip

We will use the standard $L^{2}(0,\,\ell)$ space; the scalar product and the norm are denoted by
$$
\langle\varphi,\,\psi\rangle_{L^{2}(0,\,\ell)} =
\int_{0}^{\ell}\varphi\,\overline{\psi}\,dx \text{ and }
\|\psi\|_{L^{2}(0,\,\ell)}^{2} = \int_{0}^{\ell}|\psi|^{2}dx.
$$
We introduce the Hilbert space 
\begin{eqnarray}
\label{301}\mathcal{H} & = &
H_{0}^{1}(0,\,\ell) \times H_{0}^{1}(0,\,\ell) \times H_{0}^{2}(0,\,\ell) \times L^{2}(0,\,\ell) \times L^{2}(0,\,\ell) \times L^{2}(0,\,\ell)   \nonumber \\
& &  \qquad  \times L^{2}(\mathbb{R},L^{2}(0,\,\ell)) \times L^{2}(\mathbb{R},L^{2}(0,\,\ell)) \times L^{2}(\mathbb{R},L^{2}(0,\,\ell))
\end{eqnarray}
equipped with the following  inner product:
\begin{align}
\langle {\mathcal U},\,\tilde{{\mathcal U}}\rangle_{{\mathcal H}} = &
\int_{0}^{\ell}U\overline{\tilde{U}}dx +
\int_{0}^{\ell}V\overline{\tilde{V}}dx + \int_{0}^{\ell}W\overline{\tilde{W}}dx \nonumber \\
& +\ \vartheta\int_{0}^{\ell}u_{x}\overline{\tilde{u}}_{x}dx
+ \chi\int_{0}^{\ell}v_{x}\overline{\tilde{v}}_{x}dx +  \zeta\int_{0}^{\ell}w_{xx}\overline{\tilde{w}}_{xx}dx \nonumber  \\
& +\ k\int_{0}^{\ell}(-u + v + \gamma\,w_{x})(-\overline{\tilde{u}}
+ \overline{\tilde{v}} + \gamma\overline{\tilde{w}}_{x})dx  \nonumber \\
\label{302}& +\ \mathfrak{C}
\int_{\mathbb{R}}\langle\varphi,\,\tilde{\varphi}\rangle_{L^{2}(0,\,\ell)}d\xi  +
\mathfrak{C} \int_{\mathbb{R}}\langle\phi,\,\tilde{\phi}\rangle_{L^{2}(0,\,\ell)}d\xi  + \mathfrak{C} \int_{\mathbb{R}}\langle \psi,\,\tilde{\psi}\rangle_{L^{2}(0,\,\ell)}d\xi ,
\end{align}
where ${\mathcal U}=(u,\,v,\,w,\,U,\,V,\,W,\,\varphi,\,\phi,\,\psi)^{T}$ and $\ \widetilde{{\mathcal U}}=(\tilde{u},\,\tilde{v},\,\tilde{w},\,\tilde{U},\,\tilde{V},\,\tilde{W},\,\tilde{\varphi},\,\tilde{\phi},\,\widetilde{\psi})^{T}.$ 
We  wish to transform the initial boundary value problem
\eqref{213} into an abstract Cauchy problem in the Hilbert space
$\mathcal{H}.$ For this we introduce the functions $u_{\bf t} = U,$ $v_{\bf t} = V,$ $w_{\bf t} = W,$ and we rewrite the system \eqref{213} as the following initial value
problem:
\begin{eqnarray}
\label{303}\frac{d{\mathcal U}}{d{\bf t}} ({\bf t}) = {\mathcal A}{\mathcal U}({\bf t}),\quad {\mathcal U}(0) = {\mathcal U}_{0},\quad
\forall\; {\bf t} > 0.
\end{eqnarray}
Here ${\mathcal U}$ is as above, 
${\mathcal U}_{0}=(u_{0},\,v_{0},\,w_{0},\,u_{1},\,v_{1},\,w_{0},\,0,\,0,\,0)^{T},
$ and the operator $\,\mathcal{A}:\mathcal{D}(\mathcal{A})\subset
\mathcal{H}\rightarrow \mathcal{H}$ is given by the formula
\begin{equation}
\label{304}\mathcal{A}\left(
\begin{array}{c}
u \\
\\
v\\
\\
w \\
\\
U \\
\\
V\\
\\
W \\
\\
\varphi \\
\\
\phi \\
\\
\psi
\end{array}
\right) =\left(
\begin{array}{c}
U \\
\\
V \\
\\
W \\
\\
\vartheta u_{xx} + k(-u + v + \gamma w_{x}) - \mathfrak{C}\int_{\mathbb{R}}\mu(\xi)\varphi(x,\xi)d\xi \\
\\
\chi v_{xx} - k(-u + v + \gamma w_{x}) - \mathfrak{C}\int_{\mathbb{R}}\mu(\xi)\phi(x,\xi)d\xi 
\\
\\
-\zeta w_{xxxx} + k(-u + v + \gamma w_{x})_{x} - \mathfrak{C}\int_{\mathbb{R}}\mu(\xi)\psi(x,\xi)d\xi \\
\\
-(\xi^{2} + \eta)\varphi(x,\xi) + \mu(\xi)U(x) \\
\\
-(\xi^{2} + \eta)\phi(x,\xi) + \mu(\xi)V(x) \\
\\
-(\xi^{2} + \eta)\psi(x,\xi) + \mu(\xi)W(x)
\end{array}
\right).
\end{equation}
with the domain
\begin{align*}
\mathcal{D}({\mathcal A}) = &
\left\{{\mathcal U}=(u,\,v,\,w,\,U,\,V,\,W,\,\varphi,\,\phi,\,\psi)\in {\mathcal H}:\
 U,\,V\in H_{0}^{1}(0,\,\ell), W\in H_{0}^{2}(0,\,\ell),\right. \\
\medskip
& \quad \vartheta u_{xx} + k(-u +v + \gamma w_{x}) - \mathfrak{C}\int_{\mathbb{R}}\mu(\xi)\varphi(x,\xi) d\xi \in
L^{2}(0,\,\ell), \\ 
\medskip
& \quad \chi v_{xx} - k(-u +v + \gamma w_{x}) - \mathfrak{C}\int_{\mathbb{R}}\mu(\xi)\phi(x,\xi)d\xi \in
L^{2}(0,\,\ell), \\
\medskip
& \quad - \zeta w_{xxxx} + k \gamma (-u +v + \gamma w_{x})_{x} - \mathfrak{C}\int_{\mathbb{R}}\mu(\xi)\psi(x,\xi)d\xi \in
L^{2}(0,\,\ell), \\ 
\medskip
& \quad  -(\xi^{2} + \eta)\varphi(x,\,\xi) + \mu(\xi)U(x)\in L^{2}({\mathbb{R}},L^2(0,\ell))  \\
\medskip
& \quad  -(\xi^{2} + \eta)\phi(x,\xi) + \mu(\xi)V(x)\in L^{2}({\mathbb{R}},L^2(0,\ell)) \\
\medskip
& \quad  -(\xi^{2} + \eta)\psi(x,\xi) + \mu(\xi)W(x) \in L^{2}({\mathbb{R}},L^2(0,\ell)) \\
\medskip
& \quad \ \quad \xi \varphi(x,\,\xi),\ \xi \phi(x,\,\xi),\ \xi \psi(x,\xi) \in
L^{2}({\mathbb{R}},L^2(0,\ell))\}.
\end{align*}
Let the operator $\, A :\mathcal{D}(A)\subset
H\rightarrow H$ which is given by
\begin{equation}
\label{304b} A \left(
\begin{array}{c}
u \\
v \\
w \\
\end{array}
\right) =\left(
\begin{array}{c}
- \vartheta u_{xx} - k(-u + v + \gamma w_{x}) 
\\
- \chi v_{xx} + k(-u + v + \gamma w_{x})  
\\
\zeta w_{xxxx} - k (-u + v + \gamma w_{x})_{x}  
\end{array}
\right),
\end{equation}
with the domain
$\mathcal{D}(A) =  \left[H^2(0,\ell) \cap H^1_0(0,\ell)\right]  \times  \left[H^2(0,\ell) \cap H^1_0(0,\ell)\right] \times \left[H^4(0,\ell) \cap H_{0}^{2}(0,\,\ell) \right],$  

$H = (L^2(0,\ell))^3$.

\medskip

We note that the operator $A$ is self-adjoint and strictly positive.

\medskip

So, the problem (\ref{103}) can be rewritten as following, as in \cite{kais,kaisbis}:

\begin{equation}
\label{abstract}
\left\{
\begin{array}{ll}
 Z_{\bf tt}({\bf t}) + AZ({\bf t}) + BB^* \partial_{\bf t}^{\alpha,\eta}Z ({\bf t})= 0, \ {\bf t} > 0 \\
Z(0) = Z_{0}, Z_{\bf t}(0) =  Z_{1},
\end{array}
\right.
\end{equation}

where $B = B^*=I_{H}, Z({\bf t}) = (u({\bf t}),\,v({\bf t}),\,w({\bf t}))^{T}$ and 
$Z_{0} = (u_{0},\,v_{0},\,w_{0})^{T},\quad Z_{1} = (u_{1},\,v_{1},\,w_{1})^{T}.$

\medskip

We define 

\begin{eqnarray}
\label{301bis}\mathcal{H}_{0} & = &
H_{0}^{1}(0,\,\ell) \times H_{0}^{1}(0,\,\ell) \times H_{0}^{2}(0,\,\ell) \times L^{2}(0,\,\ell) \times L^{2}(0,\,\ell) \times L^{2}(0,\,\ell)    
\end{eqnarray}
equipped with the inner product given by
\begin{align}
\langle {\mathcal U},\,\tilde{{\mathcal U}}\rangle_{{\mathcal H}_{0}} = &
\int_{0}^{\ell}U\overline{\tilde{U}}dx +
\int_{0}^{\ell}V\overline{\tilde{V}}dx + \int_{0}^{\ell}W\overline{\tilde{W}}dx \nonumber \\
& +\ \vartheta\int_{0}^{\ell}u_{x}\overline{\tilde{u}}_{x}dx
+ \chi\int_{0}^{\ell}v_{x}\overline{\tilde{v}}_{x}dx +  \zeta\int_{0}^{\ell}w_{xx}\overline{\tilde{w}}_{xx}dx \nonumber  \\
& +\ k\int_{0}^{\ell}(-u + v + \gamma\,w_{x})(-\overline{\tilde{u}}
+ \overline{\tilde{v}} + \gamma\overline{\tilde{w}}_{x})dx,
\label{302bis}
\end{align}
where ${\mathcal U}=(u,\,v,\,w,\,U,\,V,\,W)^{T}$ and $\ \tilde{{\mathcal U}}=(\tilde{u},\,\tilde{v},\,\tilde{w},\,\tilde{U},\,\tilde{V},\,\tilde{W})^{T}.$ 

\medskip

Then, the operator 
$\,\mathcal{A}_{0}:\mathcal{D}(\mathcal{A}_{0})\subset
\mathcal{H}_{0}\rightarrow \mathcal{H}_{0}$ given by
\begin{equation}
\label{304bis}\mathcal{A}_{0} \left(
\begin{array}{c}
u \\
\\
v\\
\\
w \\
\\
U \\
\\
V\\
\\
W \\
\end{array}
\right) =\left(
\begin{array}{c}
U \\
\\
V \\
\\
W \\
\\
\vartheta u_{xx} + k(-u + v + \gamma w_{x}) - U \\
\\
\chi v_{xx} - k(-u + v + \gamma w_{x}) - V
\\
\\
-\zeta w_{xxxx} + k(-u + v + \gamma w_{x})_{x} - W  
\end{array}
\right).
\end{equation}
with the domain
\begin{align*}
\mathcal{D}({\mathcal A}_{0}) =  [H^{2}(0,\,\ell) \times H_{0}^{1}(0,\,\ell)]^{2} \times [H^{4}(0,\,\ell) \cap H_{0}^{2}(0,\,\ell)] \times H_{0}^{2}(0,\,\ell).
\end{align*}

We recall two propositions from \cite{kais,kaisbis}:
 
\begin{proposition}{\cite[Proposition 2.6.2]{kaisbis}}
The operator ${\mathcal A}_{0}$ generates a C$_{0}$-semigroup of contractions in the Hilbert space ${\mathcal H}_{0}.$ Moreover, the following auxiliary problem:
\begin{eqnarray}
\label{606} \left\lbrace
\begin{array}{l}
u_{\bf tt}({\bf t}) - \vartheta u_{xx}({\bf t}) - k(-u({\bf t}) + v({\bf t}) + \gamma w_{x}({\bf t})) + u_{\bf t}({\bf t}) = 0,   \\
v_{\bf tt}({\bf t}) - \chi v_{xx}({\bf t}) + k(-u({\bf t}) + v({\bf t}) + \gamma w_{x}({\bf t})) 
+ v_{\bf t}({\bf t}) = 0,   \\
w_{\bf tt}({\bf t}) + \zeta w_{xxxx}({\bf t}) - k\gamma(-u({\bf t}) + v({\bf t}) + \gamma w_{x}({\bf t}))_{x} + w_{\bf t}({\bf t}) = 0,  \\
u(0,\,{\bf t}) = u(\ell,\,{\bf t}) = 0,\ v(0,\,{\bf t}) = v(\ell,\,{\bf t}) = 0, \\ 
w(0,\,{\bf t}) = w(\ell,\,{\bf t}) = w_{x}(0,\,{\bf t}) = w_{x}(\ell,\,{\bf t}) = 0,
\\
u(x,\,0)=u_{0}(x),\quad u_{\bf t}(x,\,0) = u_{1}(x),   \\
v(x,\,0)=v_{0}(x),\quad v_{\bf t}(x,\,0) = v_{1}(x),   \\
w(x,\,0) = w_{0}(x),\quad w_{\bf t}(x,\,0) = w_{1}(x),
\end{array}
\right.
\end{eqnarray}
admits a unique solution $(u(x,\,{\bf t}),\,v(x,\,{\bf t}),\,w(x,\,{\bf t}))$ such that if $(u_{0},\,v_{0},\,w_{0},\,u_{1},\,v_{1},\,w_{1})\in {\mathcal D}({\mathcal A}_{0})$, then the solution $(u(x,\,{\bf t}),\,v(x,\,{\bf t}),\,w(x,\,{\bf t}))$ of \eqref{606} verifies the following regularity property:
\begin{align*}
{\mathcal U}=(u,\,v,\,w,\,u_{\bf t},\,v_{\bf t},\,w_{\bf t})\in C([0,\,+\infty):\,{\mathcal D}({\mathcal A}_{0}))\cap C^{1}([0,\,+\infty),\,{\mathcal H}_{0}),
\end{align*}
and when $(u_{0},\,v_{0},\,w_{0},\,u_{1},\,v_{1},\,w_{1})\in {\mathcal H}_{0}$, then 
\begin{align*}
{\mathcal U}=(u,\,v,\,w,\,u_{\bf t},\,v_{\bf t},\,w_{\bf t})\in C([0,\,+\infty),\,{\mathcal H}_{0}).
\end{align*}
\end{proposition}
The energy of the solution of the system \eqref{606} is defined as follows:
\begin{align}
E_{0}({\bf t}) = &\ \frac{1}{2}\left[\|u_{\bf t}({\bf t})\|_{L^{2}(0,\,\ell)}^{2} +
\|v_{\bf t}({\bf t})\|_{L^{2}(0,\,\ell)}^{2} + \|w_{\bf t}({\bf t})\|_{L^{2}(0,\,\ell)}^{2} 
\right. \nonumber  \\
&\qquad +\ \vartheta\|u({\bf t})\|_{H_{0}^{1}(0,\,\ell)}^{2} +
\chi\|v({\bf t})\|_{H_{0}^{1}(0,\,\ell)}^{2} + \zeta\|w({\bf t})\|_{H_{0}^{2}(0,\,\ell)}^{2} \nonumber \\
\label{607}&\qquad\left. +\ k\|-u({\bf t}) + v({\bf t}) + \gamma w_{x}({\bf t})\|_{L^{2}(0,\,\ell)}^{2}\right],
\end{align}
and it is decreasing function of the time ${\bf t}.$ In particular, we have
\begin{align}
\label{608}\frac{dE_{0}({\bf t})}{d{\bf t}} = -\|u_{\bf t}({\bf t})\|_{L^{2}(0,\,\ell)}^{2} - \|v_{\bf t}({\bf t})\|_{L^{2}(0,\,\ell)}^{2} - \|w_{\bf t}({\bf t})\|_{L^{2}(0,\,\ell)}^{2}.
\end{align}

\begin{proposition}{\cite[Theorem 2.3.1]{kaisbis}}
\label{proposition301}
The operator ${\mathcal A}$ generates a C$_{0}$-semigroup of contractions in the Hilbert space ${\mathcal H}.$ Hence the problem \eqref{606} admits a unique solution 
$$
(u(x,\,{\bf t}),\,v(x,\,{\bf t}),\,w(x,\,{\bf t}),\varphi(x,\,{\bf t},\,\xi),\,\phi(x,\,{\bf t},\,\xi),\,\psi(x,\,{\bf t},\,\xi))
$$ 
such that if $(u_{0},\,v_{0},\,w_{0},\,u_{1},\,v_{1},\,w_{1},\,0,\,0,\,0)\in {\mathcal D}({\mathcal A})$, then the solution $(u,\,v,\,w,\varphi,\,\phi,\,\psi)$ of \eqref{606} satisfies the following regularity property:
\begin{align*}
{\mathcal U}=(u,\,v,\,w,\,u_{\bf t},\,v_{\bf t},\,w_{\bf t},\,\varphi,\,\phi,\,\psi)\in C([0,\,+\infty),\,{\mathcal D}({\mathcal A}))\cap C^{1}([0,\,+\infty),\,{\mathcal H}),
\end{align*}
and when $(u_{0},\,v_{0},\,w_{0},\,u_{1},\,v_{1},\,w_{1},\,0,\,0,\,0)\in {\mathcal H}$, then 
\begin{align*}
{\mathcal U}=(u,\,v,\,w,\,u_{\bf t},\,v_{\bf t},\,w_{\bf t},\,\varphi,\,\phi,\,\psi)\in C([0,\,+\infty),\,{\mathcal H}).
\end{align*}
\end{proposition}
The energy of the system \eqref{606}, defined by
\begin{align}
E({\bf t}) := \frac{1}{2} \left\|(u,\,v,\,w,\,u_{\bf t},\,v_{\bf t},\,w_{\bf t},\,\varphi,\,\phi,\,\psi)\right\|_{\mathcal{H}}^2
\label{607bis},
\end{align}
is decreasing in  time ${\bf t}.$ In particular, we have
\begin{align}
\label{608bis}\frac{dE({\bf t})}{d{\bf t}} = & - \mathfrak{C}\int_\mathbb{R}(|\xi|^{2} + \eta)\|\varphi({\bf t,\,\xi})\|_{L^{2}(0,\,\ell)}^{2}d\xi - \mathfrak{C}\int_{\mathbb{R}}(|\xi|^2 + \eta)\|\phi({\bf t,\,\xi})\|_{L^{2}(0,\,\ell)}^{2}d\xi \nonumber\\
& - \mathfrak{C}\int_\mathbb{R}(|\xi|^{2} + \eta)\|\psi({\bf t,\,\xi})\|_{L^{2}(0,\,\ell)}^{2}d\xi.
\end{align}
 
\section{Non-uniform stabilization} \label{4}
\setcounter{equation}{0}

\subsection{Polynomial stability for $\eta =0.$}
Assume that $\eta = 0,$ then the operator ${\mathcal A}$ is not onto and consequently $0\notin \varrho({\mathcal A})$ the resolvent set of ${\mathcal A}.$

Since the embedding $H_{0}^{1}(0,\,\ell) \times H_{0}^{1}(0,\,\ell) \times H_{0}^{2}(0,\,\ell) \subset L^{2}(0,\,\ell) \times L^{2}(0,\,\ell) \times L^{2}(0,\,\ell)$ is compact and the only solution of the following problem:

\begin{eqnarray}
\left\lbrace
\begin{array}{l}
u_{\bf tt} - \vartheta u_{xx} - k(-u + v + \gamma w_{x}) = 0,\quad  (0,\,\ell) \times(0,\,+\infty),\\
v_{\bf tt} - \chi v_{xx} + k(-u + v + \gamma w_{x})  = 0,   \quad (0,\,\ell) \times(0,\,+\infty),\\
w_{\bf tt} + \zeta w_{xxxx} - k (-u + v + \gamma w_{x})_{x}  = 0,\quad (0,\,\ell) \times (0,\,+\infty),\\
u_{\bf t}(x,{\bf t}) = 0, v_{\bf t}(x,{\bf t}) =0, w_{\bf t}(x,{\bf t}) = 0, \ (0,\,\ell) \times(0,\,+\infty),\\
u(0,\,{\bf t}) = u(\ell,\,{\bf t}) = 0,\ v(0,\,{\bf t}) = v(\ell,\,{\bf t}) = 0,\quad  (0,\,+\infty), \\ 
w(0,\,{\bf t}) = w(\ell,\,{\bf t}) = w_{x}(0,\,{\bf t}) = w_{x}(\ell,\,{\bf t}) = 0, (0,\,+\infty),
\\
u(x,\,0)=u_{0}(x),\quad u_{\bf t}(x,\,0) = u_{1}(x),\quad (0,\,\ell),  \\
v(x,\,0)=v_{0}(x),\quad v_{\bf t}(x,\,0) = v_{1}(x),\quad (0,\,\ell),  \\
w(x,\,0) = w_{0}(x),\quad w_{\bf t}(x,\,0) = w_{1}(x), \quad (0,\,\ell),
\end{array}
\right.
\end{eqnarray}
is the trivial one.

\medskip

Then according to \cite[Sect. 2.4]{kaisbis} the semigroup $(e^{{\bf t}\mathcal{A}})_{t \geq 0}$ is strongly stable, i.e.,
$$\|e^{{\bf t}{\mathcal A}}{\mathcal U}_{0}\|_{{\mathcal H}} \longrightarrow 0\quad {\rm as}\quad {\bf t}\rightarrow +\infty$$
for all initial data ${\mathcal U}_{0}\in {\mathcal H}$.

To determine the rate of stability in this case, we need the following result \cite[Theorem 8.4]{BCT}: 

\begin{theorem}[\cite{BCT}] \label{TBC}
Let $\mathcal{T}(t)$ be a bounded C$_0$-semigroup on a Hilbert space $X$ with generator $H.$ Assume that $\sigma (H) \cap i \mathbb{R} = \left\{ 0\right\}$ and that there exist $\beta \geq 1, \gamma > 0$ such that 
$$
\left\| (i\lambda I - H)^{-1}\right\|_{\mathcal{L}(X)} \leq \left\{
\begin{array}{ll}
O(|\lambda|^{-\beta}), \, \lambda \rightarrow 0, \\
O(|\lambda|^{\gamma}), \, \lambda \rightarrow \infty.
\end{array}
\right.
$$
Then, there exists  a constant $C > 0$ such that  
$$
\left\|\mathcal{T}(t) z \right\|_{X} \leq 
\frac{C}{{\bf t}^{\frac{1}{\max{(\beta,\gamma)}}}} \left\|z\right\|_{{\mathcal D}(H)} , \; \forall \, t > 0, z \in {\mathcal D}(H) \cap R(H),
$$
where ${\mathcal D}(H)$ is the domain of $H$ and $R(H)$ is the range of $H.$
 \end{theorem}

Based on Theorem \ref{TBC}, a simple adaptation of the proofs of \cite[Theorems 5.8 and 5.9]{abbes} (see also \cite{kais,kaisbis}) lead to the following stability result: 

 \begin{proposition}
If $\eta=0$, then there exists a $C > 0$ such that 
\begin{equation}
\label{eta}
\left\|e^{{\bf t}\mathcal{A}} U_0\right\|_{\mathcal{H}}  \leq 
\frac{C}{\sqrt{{\bf t}}} \left\|U_0\right\|_{\mathcal{D}(\mathcal{A})}, \, \forall \, {\bf t} >  0, \, U_0 \in {\mathcal D}(\mathcal{A}) \cap R(\mathcal{A}),
\end{equation}
where $R(\mathcal{A})$ is the range of $\mathcal{A}.$
 \end{proposition}

\subsection{Polynomial stability for $\eta > 0$}
In this case there exists $\delta,\ C > 0$, such that the auxiliary dissipative operator satisfies the following relation: 
$$
\left\|e^{{\bf t}\mathcal{A}_0}\right\|_{\mathcal{L}(\mathcal{H}_{0})} \leq C\,e^{-\delta {\bf t}},\quad \forall \ {\bf t} \geq 0.
$$
Applying \cite[Corollary 2.6.1]{kaisbis} hence we obtain the following polynomial decay result for the systems \eqref{103} and \eqref{213}:

\begin{proposition} \label{corollary905}
 If $\eta > 0$, then the semigroup $(e^{{\bf t}{\mathcal A}})_{t \geq 0}$ is polynomially stable, namely there exists a constant $C>0$ such that
\begin{align*}
& \|e^{{\bf t}{\mathcal A}}(u_{0},\,v_{0},\,w_{0},\,u_{1},\,v_{1},\,w_{1},\,\varphi_{0},\,\phi_{0},\,\psi_{0})\|_{{\mathcal H}} \\
& \qquad\qquad \leq \frac{C}{(1 + {\bf t})^{\frac{1}{1 - \alpha}}}\|(u_{0},\,v_{0},\,w_{0},\,u_{1},\,v_{1},\,w_{1},\,\varphi_{0},\,\phi_{0},\,\psi_{0})\|_{{\mathcal D}({\mathcal A})},\quad \forall\;{\bf t}\geq 0,
\end{align*}
for all initial data $(u_{0},\,v_{0},\,w_{0},\,u_{1},\,v_{1},\,w_{1},\,\varphi_{0},\,\phi_{0},\,\psi_{0}) \in {\mathcal D}({\mathcal A}).$ In particular, the energy of the strong solution of \eqref{103} and \eqref{213} satisfies the following estimate:
\begin{align*}
E({\bf t}) \leq \frac{C}{(1 + {\bf t})^{\frac{2}{1 - \alpha}}}\|(u_{0},\,v_{0},\,w_{0},\,u_{1},\,v_{1},\,w_{1},\,0,\,0,\,0)\|_{{\mathcal D}({\mathcal A})}^{2}, \quad \forall  \, {\bf t} >0.
\end{align*} 
\end{proposition} 

\section{Numerical study} \label{5}
\setcounter{equation}{0}


In this section we will verify numerically the polynomial rate of decay obtained in the previous 
section. 

\subsection{Linear equations of Motion}

First, we approximate  the longitudinal and transverse displacement vector $[u,v,w]^\top$
in space using an energy-conservative finite difference method.
For $J \in \mathbb{N}$ and $\delta x=\ell/J$,
we define
$x_{j}$, with $j=1,\ldots,J$, as a uniform discretization of the interval
$(0,\ell)$,
obtaining a vector
$\mathbf{U}=[\mathbf{u}(t), \mathbf{v}(t), \mathbf{w}(t)]^\top$
approximation of $[u,v,w]^\top$
in $\mathbb{R}^{3J}$.
We obtain 
the linear equation of motion
\begin{equation}
	\label{LEM}
	\mathbf{M}\ddot{\mathbf{U}}(t)+
	\mathbf{K}{\mathbf{U}}(t)+
	\mathbf{C}\overset{\alpha,\eta}{\mathbf{U}}(t)
	=0
\end{equation}
where $\overset{\alpha,\eta}{\mathbf{U}}(t)= D^{\alpha,\eta}{\mathbf{U}}(t)$
is the generalized Caputo fractional derivative defined in \eqref{202}.
In addition, and according to a discretization consistent with
\eqref{103}, we choose
$\mathbf{M}= \mathbf{C} = \mathbf{I}_{3J}$ (identity matrix).
The stiffness matrix is given by $\mathbf{K}=\mathbf{K}_{\mathrm{stress}} +  \mathbf{K}_{\mathrm{coupling}}$, where
$$
\mathbf{K}_{\mathrm{stress}} = 
\begin{pmatrix}
-\vartheta \mathbf{D}^2 & & \\ & -\chi\mathbf{D}^2 & \\ & & \zeta\mathbf{D}^4 
\end{pmatrix},
\qquad
\mathbf{K}_{\mathrm{coupling}} = 
k\begin{pmatrix}
 \mathbf{I}_J & - \mathbf{I}_J  & \gamma \mathbf{R} \\  -\mathbf{I}_J &  \mathbf{I}_J  & -\gamma \mathbf{R} & \\ 
 \gamma \mathbf{R}^T & -\gamma \mathbf{R}^T & -\gamma^2  \mathbf{D}^2.
\end{pmatrix}
$$
Here, $\mathbf{D}^2$ and $\mathbf{D}^4$ are the finite difference matrix approximation of
$\partial_{xx}$ and $\partial_{xxxx}$ respectively with the boundary conditions \eqref{103}$_{7,8}$,
and $\mathbf{R}$ is the upper triangular matrix, obtained from the Cholesky decomposition
$\mathbf{D}^2=- \mathbf{R}^T\mathbf{R}$.
Using uniform finite differences we have
\begin{eqnarray*}
   \left(\mathbf{D}^2\mathbf{u}\right)_j
   &=&
 \dfrac{\mathbf{u}_{j-1}
 -2\mathbf{u}_{j}
  +\mathbf{u}_{j+1}
 }
 {\delta x^2},
 \\
   \left(\mathbf{D}^4\mathbf{w}\right)_j&=&
 \dfrac{\mathbf{w}_{j-2}
 -4\mathbf{w}_{j-1}
 +6\mathbf{w}_{j}
 -4\mathbf{w}_{j+1}
 +\mathbf{w}_{j+2}
 }
 {\delta x^4},\qquad J=1,\ldots,J.
 \label{D4}
\end{eqnarray*}
According to the boundary conditions \eqref{103}$_{7,8}$, we seek the solution of 
Finite Difference in
$$
X_J:=\left\{ 
\mathbf{U}=[\mathbf{u}(t), \mathbf{v}(t), \mathbf{w}(t)]^\top \in \mathbb{R}^{3J} \right\},
$$
considering
$
\mathbf{u}_0=\mathbf{u}_{J+1}=\mathbf{v}_0=\mathbf{v}_{J+1}=
\mathbf{w}_0=\mathbf{w}_{J+1}=\mathbf{v}_{-1}=\mathbf{v}_{J+2}=0
$.

\subsection{Time discretization}

In order to preserve the energy with a second order scheme in time, we choose a $\beta$-Newmark scheme for $w$.
The method  consists of updating
the displacement, velocity and acceleration vectors
 at the current time $t^n=n\delta t$ to the time $t^{n+1} = (n+1)\delta t$, a small time interval
 $\delta t$
later.
The Newmark algorithm \cite{Newmark}
is based on a set of two relations expressing the forward displacement
$\mathbf{U}^{n+1}$
 and velocity
$\dot{\mathbf{U}}^{n+1}$  in terms of their
current values and the forward and current values of the acceleration:
\begin{eqnarray}
\dot{\mathbf{U}}^{n+1} &=& \dot{\mathbf{U}}^{n} + (1-\widetilde{\gamma})\delta t\,\ddot{\mathbf{U}}^{n} + \widetilde{\gamma}\delta t\,\ddot{\mathbf{U}}^{n+1},\label{702}\\
{\mathbf{U}}^{n+1} &=& {\mathbf{U}}^{n} + \delta t \dot{\mathbf{U}}^{n} + \left(\frac{1}{2}-\widetilde{\beta}\right)\delta t^ 2\,
\ddot{\mathbf{U}}^{n} + \widetilde{\beta}\delta t^2\,\ddot{\mathbf{U}}^{n+1},\label{703}
\end{eqnarray}
where $\widetilde{\beta}$ and $\widetilde{\gamma}$ are parameters of the methods that will be fixed later.
Replacing \eqref{702}--\eqref{703} in the equation of motion \eqref{LEM}, we obtain
\begin{equation}
      \left(
    \mathbf{M}
    +  \widetilde{\beta}\delta t^2\,
    \mathbf{K}
\right)
\ddot{\mathbf{U}}^{n+1}
+ 
	\mathbf{C} \halfscript{\overset{\alpha,\eta}{\mathbf{U}}}{n+1}
=
 -
  \mathbf{K}
\left( \mathbf{U}^n+\delta t \dot{\mathbf{U}}^n + \left(\dfrac{1}{2}- \widetilde{\beta}\right) \delta t^2 \ddot{\mathbf{U}}^n\right).
 \label{NewmarkW}
\end{equation}

\subsection{Approximation of fractional derivatives}
At this state we have two possibilities: (1) approximate the fractional derivative using a classical finite difference method such as the truncation of the Grünwald–-Letnikov derivative \cite{Ortigueira}, which in turn is an equivalent form of the Riemann--Liouville derivative;
or, (2) propose a conservative numerical scheme for the Mbodje augmented  model \eqref{213}.

\subsubsection{Approximation using a truncation of the Grünwald-–Letnikov derivative.}
In the case  of the classical Grünwald–Letnikov derivative 
we consider a fractional trapezoidal formula \cite{li}, in order to approximate
\begin{equation}
\left[\prescript{}{C}D^{\alpha,\eta}_{0,t}f\right]_{t=t_n} \approx \sum_{k=0}^n a_{k,n} e^{-\eta\delta t (n-k)}f'(t_k),
\label{caputoapprox}
\end{equation}
where 
\begin{equation}
a_{k,n}=\dfrac{\delta t^{1-\alpha}}{\Gamma(3-\alpha)}
\begin{cases}
(n-1)^{2-\alpha}-(n-2+\alpha)n^{1-\alpha},& k=0,\\
(n-k+1)^{2-\alpha}+(n-1-k)^{2-\alpha}-2(n-k)^{2-\alpha},& 1\leqslant k \leqslant n-1,\\
1,& k=n.
\end{cases}
\label{akn}
\end{equation}
Replacing the term $\halfscript{\overset{\alpha,\eta}{\mathbf{U}}}{n+1}$ by the numerical approximation $\left[\prescript{}{C}D^{\alpha ,\eta}_{0,t}f\right]_{t=t_n}$ given by \eqref{caputoapprox}--\eqref{akn}, then the equation \eqref{NewmarkW} becomes
\begin{eqnarray}
\lefteqn{
      \left(
    \mathbf{M}
    +  \widetilde{\beta}\delta t^2\,
    \mathbf{K}
    +
     \dfrac{\widetilde{\gamma}\delta t^{2-\alpha}}{\Gamma(3-\alpha)}\,
    \mathbf{C}
\right)
\ddot{\mathbf{U}}^{n+1}
=
 -
  \mathbf{K}
\left( \mathbf{U}^n+\delta t \dot{\mathbf{U}}^n + \left(\dfrac{1}{2}- \widetilde{\beta}\right) \delta t^2 \ddot{\mathbf{U}}^n\right)}
\nonumber
\\
&& \qquad\qquad\qquad\qquad
-    \dfrac{\delta t^{1-\alpha}}{\Gamma(3-\alpha)}\,
    \mathbf{C}
\left(  \dot{\mathbf{U}}^n + \left(1- \widetilde{\gamma}\right) \delta t \ddot{\mathbf{U}}^n
+\sum_{k=0}^n a_{k,n+1} \dot{\mathbf{U}}^k
\right).
 \label{Unmas1}
\end{eqnarray}

\subsubsection{Approximation using the Mbodje augmented  model}
On the other hand, in  case of  the Mbodje augmented  model \eqref{213}
we propose the following  conservative numerical scheme:
\begin{equation}
    \begin{cases}
      \left(
    \mathbf{M}
    +  \widetilde{\beta}\delta t^2\,
    \mathbf{K}
\right)
\ddot{\mathbf{U}}^{n+1}
+ 
	\mathfrak{C}  \displaystyle\sum_{\ell=1}^M\mu_\ell \mathbf{\Phi}^{n+1}_\ell
=
 -
  \mathbf{K}
\left( \mathbf{U}^n+\delta t \dot{\mathbf{U}}^n + \left(\dfrac{1}{2}- \widetilde{\beta}\right) \delta t^2 \ddot{\mathbf{U}}^n\right),\\
\mathbf{\Phi}^{n+1}_\ell = \mathbf{\Phi}^{n}_\ell
- \delta t \left(\xi_\ell^2 +\eta\right) \mathbf{\Phi}^{n+\frac{1}{2}}_\ell
+ \delta t \mu_\ell  \dot{\mathbf{U}}^{n+\frac{1}{2}},
\end{cases}
 \label{ConservNumerical}
\end{equation}
where 
$\mathbf{\Phi}^{n+\frac{1}{2}}_\ell
=\dfrac{\mathbf{\Phi}^{n}_\ell+\mathbf{\Phi}^{n+1}_\ell}{2}
$, with
$\mathbf{\Phi}^n_\ell$ denoting the
approximation of $[\phi,\varphi,\psi]^\top
\in \mathbb{R}^{3J}$ evaluated in $\xi_\ell:=\ell\delta \xi$
and $t=n\delta t$,
for $\ell=1,\ldots,M$, $n=1,\ldots,N$
and a fixed $\delta\xi>0$.
Using  \eqref{702} again and replacing in 
\eqref{ConservNumerical}$_2$,
we can rewrite \eqref{ConservNumerical}
in the following more explicit and computable way:
\begin{equation}
    \begin{cases}
      \left(
    \mathbf{M}
    +  \widetilde{\gamma}\delta t\,
    \mathbf{C}_{augm}
    +  \widetilde{\beta}\delta t^2\,
    \mathbf{K}
\right)
\ddot{\mathbf{U}}^{n+1}
=
-
	\mathfrak{C}  \displaystyle\sum_{\ell=1}^M\widetilde{\mu}_\ell \mathbf{\Phi}^{n}_\ell
 -\mathbf{C}_{augm}\left(
\dot{\mathbf{U}}^n + \left( 1 - \widetilde{\gamma} \right)
\delta t \ddot{\mathbf{U}}^n
\right)
  \\
\qquad\qquad\qquad\qquad\qquad\qquad\qquad - 
\mathbf{K}
\left( \mathbf{U}^n+\delta t \dot{\mathbf{U}}^n + \left(\dfrac{1}{2}- \widetilde{\beta}\right) \delta t^2 \ddot{\mathbf{U}}^n\right),
\\
\mathbf{\Phi}^{n+1}_\ell = \dfrac{2-\delta t\left(\xi_\ell^2 +\eta\right)}{2+\delta t\left(\xi_\ell^2 +\eta\right)} \mathbf{\Phi}^{n}_\ell
+ 
\dfrac{2\delta t \mu_\ell }{2+\delta t\left(\xi_\ell^2 +\eta\right)}
\dot{\mathbf{U}}^{n+\frac{1}{2}},
\end{cases}
 \label{ExplicitWay}
\end{equation}
where $\mu_\ell=\xi_\ell^{(2\alpha-1)/2}$, $\widetilde{\mu}_\ell=\dfrac{2-\delta t\left(\xi_\ell^2 +\eta\right)}{2+\delta t\left(\xi_\ell^2 +\eta\right)}
\mu_\ell$ and $\mathbf{C}_{augm}=\delta t \mathfrak{C}  \left(\displaystyle\sum_{\ell=1}^M
\dfrac{2\widetilde{\mu}_\ell^2}{2+\delta t\left(\xi_\ell^2 +\eta\right)}\right)\textbf{I}_{3J} $.
Although the two approximations are consistently valid from the point of view of finite differences, we will opt for the second approximation because it conserves the energy in the non-dissipative case, and its decrease is consistent with the continuous case. It also involves the calculation of $\mathbf{\Phi}$ which is inherent to the definition of energy itself. 

\subsection{Decay of the discrete energy using the Mbodje augmented  model}

Evaluating \eqref{ConservNumerical}$_1$ in $t=t_{n+\frac{1}{2}}$, multiplying 
by $\ddot{\mathbf{U}}^{n+\frac{1}{2}}$, and summing
\eqref{ConservNumerical}$_2$ multiplied by 
$\mathfrak{C}\mathbf{\Phi}^{n+\frac{1}{2}}_\ell$, we obtain that
\begin{eqnarray}
\lefteqn{
\left[ E^{\text{Mbodje}}_\Delta \right]_n^{n+1}:=
 \left[ 
 \dfrac{1}{2} \dot{\mathbf{U}}^T  \mathbf{M} \dot{\mathbf{U}}
 +
  \dfrac{1}{2} {\mathbf{U}}^T  \mathbf{K} {\mathbf{U}}
  +
 \dfrac{\mathfrak{C}}{2}  \displaystyle\sum_{\ell=1}^M {\mu}_\ell \left|\mathbf{\Phi}_\ell\right|^2
  \right]_n^{n+1}
  }\nonumber \\
  &=&
   \Bigg[ 
   - \dfrac{\vartheta}{2} \dot{\mathbf{u}}^T  \mathbf{D}^2 \dot{\mathbf{u}}
  - \dfrac{\chi}{2} \dot{\mathbf{v}}^T  \mathbf{D}^2 \dot{\mathbf{v}}
  +\dfrac{\zeta}{2} \dot{\mathbf{w}}^T  \mathbf{D}^4 \dot{\mathbf{w}}
   +
 \dfrac{k}{2} \| - \mathbf{u} + \mathbf{v} + \gamma \mathbf{R}\mathbf{w}  \|^2
  \nonumber \\
  &&
 +\ 
 \dfrac{\mathfrak{C}}{2}  \displaystyle\sum_{\ell=1}^M {\mu}_\ell \left(\left|{\phi}_\ell\right|^2+\left|{\varphi}_\ell\right|^2+
 \left|{\psi}_\ell\right|^2\right)
  \Bigg]_n^{n+1}  \nonumber \\
  &=&
  - \ \mathfrak{C}
   \displaystyle\sum_{\ell=1}^M  \left(\xi_\ell^2 +\eta\right) \left(\left|{\phi}^{n+\frac{1}{2}}_\ell\right|^2+\left|{\varphi}^{n+\frac{1}{2}}_\ell\right|^2+
 \left|{\psi}^{n+\frac{1}{2}}_\ell\right|^2\right)
   \label{ener_Mbodje}
\end{eqnarray}
which is consistent with the estimates \eqref{607bis}--\eqref{608bis}, and constitutes a correct approximation of the energy and its decreasing behavior.

\begin{remark}
In the case of an approximation using a truncation of theGrünwald–-Letnikov derivative, we have the following estimate of the energy.
Evaluating \eqref{LEM} in $t=t_{n+\frac{1}{2}}$ and multiplying 
by $\ddot{\mathbf{U}}^{n+\frac{1}{2}}$ we obtain
\begin{eqnarray}
\lefteqn{
\left[ E^{GL}_\Delta \right]_n^{n+1}:=
 \left[ 
 \dfrac{1}{2} \dot{\mathbf{U}}^T  \mathbf{M} \dot{\mathbf{U}}
 +
  \dfrac{1}{2} {\mathbf{U}}^T  \mathbf{K} {\mathbf{U}}
  \right]_n^{n+1}
  }\nonumber \\
  &&=
   \left[ 
   - \dfrac{\vartheta}{2} \dot{\mathbf{u}}^T  \mathbf{D}^2 \dot{\mathbf{u}}
  - \dfrac{\chi}{2} \dot{\mathbf{v}}^T  \mathbf{D}^2 \dot{\mathbf{v}}
  +\dfrac{\zeta}{2} \dot{\mathbf{w}}^T  \mathbf{D}^4 \dot{\mathbf{w}}
  +
 \dfrac{k}{2} \| - \mathbf{u} + \mathbf{v} + \gamma \mathbf{R}\mathbf{w}  \|^2
  \right]_n^{n+1}  \nonumber \\
  &&=
  - \frac{\delta t^{2-\alpha}}{\Gamma(3-\alpha)} 
   \left(\dot{\mathbf{U}}^{n+\frac{1}{2}}\right)^T  \mathbf{C} \dot{\mathbf{U}}^{n+\frac{1}{2}}
    - (2^{1-\alpha}-1) \frac{\delta t^{2-\alpha}}{\Gamma(3-\alpha)} 
   \left(\dot{\mathbf{U}}^{n+\frac{1}{2}}\right)^T  \mathbf{C} \dot{\mathbf{U}}^{n}
  \nonumber \\ &&~~~
   - \sum_{k=0}^{n-1} a_{k,n+\frac{1}{2}}\left(\dot{\mathbf{U}}^{n+\frac{1}{2}}\right)^T  \mathbf{C} \dot{\mathbf{U}}^{k}
   \label{ener_est}
\end{eqnarray}
where $ a_{k,n+\frac{1}{2}} = \dfrac{a_{k,n+1}+ a_{k,n}}{2}$.
The first term on the right hand side of \eqref{ener_est} is obviously strictly negative. However, the sign of the second and third terms in \eqref{ener_est}  cannot be guaranteed. In fact, it is verified numerically that the energy discretized in this way is not necessarily strictly decreasing, unless the solution is sufficiently smooth and delta t is small enough.
\begin{figure}
    \centering
    \includegraphics[scale=0.35]{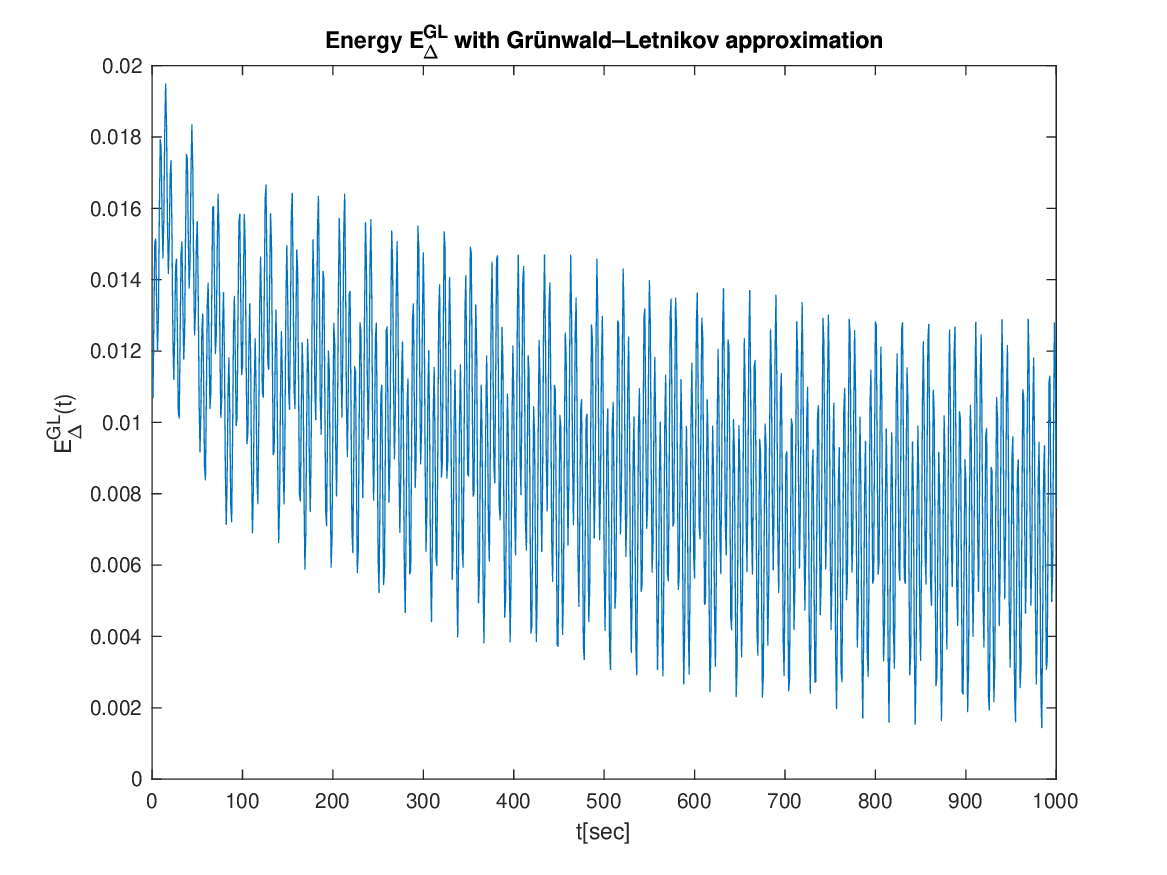}
    \includegraphics[scale=0.35]{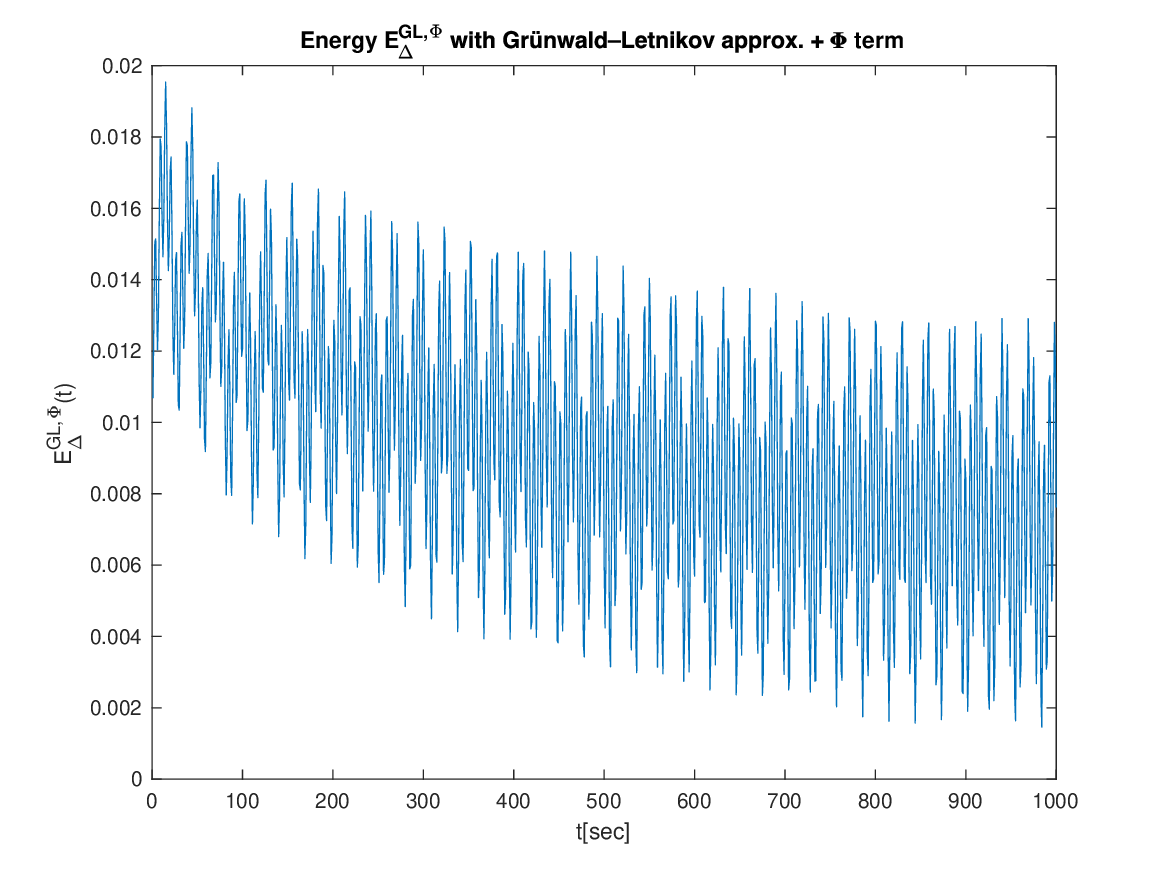}\\
    \includegraphics[scale=0.35]{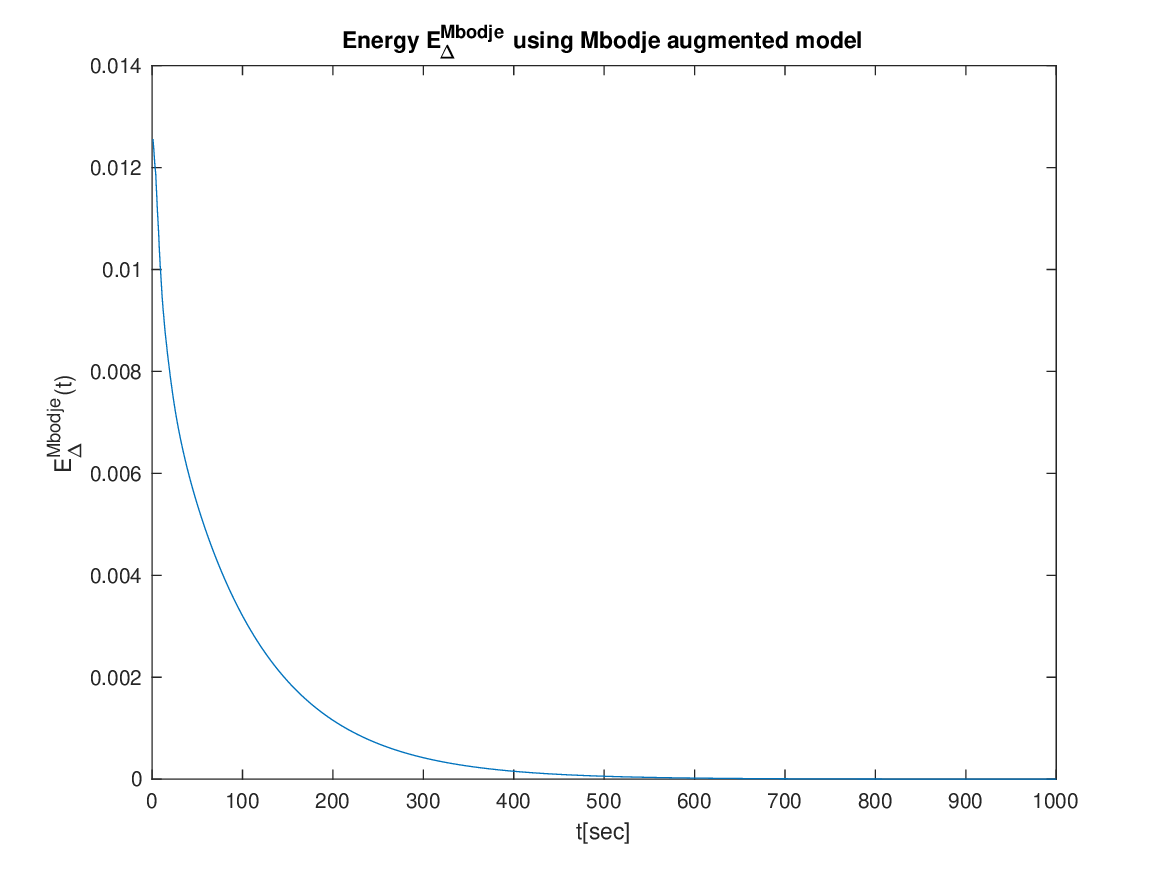}
    \caption{Comparison between the Grünwald–Letnikov derivative
    approximation and the discretized Mbodje augmented model.
    At the top left: $E^{GL}_\Delta(t)=\dfrac{1}{2} \dot{\mathbf{U}}^T  \mathbf{M} \dot{\mathbf{U}}
 +
  \dfrac{1}{2} {\mathbf{U}}^T  \mathbf{K} {\mathbf{U}}$;
  At the top right: $E^{GL,\Phi}_\Delta(t)=E^{GL}_\Delta(t)
  +
 \dfrac{\mathfrak{C}}{2}  \displaystyle\sum_{\ell=1}^M {\mu}_\ell \left|\mathbf{\Phi}_\ell\right|^2$;
 At the bottom: $E^{Mbodje}_\Delta(t)$ defined in \eqref{ener_Mbodje}.
  }
    \label{fig:enter-label}
\end{figure}
Figure 1 shows a comparison of both models. It is proved that the calculation of energy using the classical Grünwald–-Letnikov derivatives fails, either by using the derivative defined in \eqref{ener_est}, (in top left graph), as well as by artificially adding the calculation of $\mathbf{\Phi}$ as seen in the top right graph. Therefore, the only acceptable method for all our calculations should be the discretization of the augmented Mbodje model (bottom graph).
\end{remark}

\subsection{Numerical examples}

We consider here the value of the parameters $\vartheta=\chi=\zeta=k=\gamma=1$, a beam of length $\ell=1$,
and the initial conditions
\begin{equation}
    \label{CI}
\begin{cases}
\displaystyle u(x,0)=x\left(
\frac{\ell^3}{8}-\left|x-\frac{\ell}{2}\right|^3\right), & u_t(x,0)=0,\\
\displaystyle v(x,0)=0, & v_t(x,0)=0,\\
w(x,0)=x\left(x-\frac{\ell}{2}\right)\left(
\displaystyle \frac{\ell^3}{8}-\left|x-\frac{\ell}{2}\right|^3\right), & w_t(x,0)=0.
\end{cases}
\end{equation}
We observe that $u(0,t)=u(\ell,t)=0$ and $u(\cdot,0)$ is of class $\mathcal{C}^1(0,\ell)$ but its second derivative has a discontinuity in $x=\ell/2$, and in turn, $w(0,t)=w(\ell,t)=0$ with $w(\cdot,0)$ of class $\mathcal{C}^3(0,\ell)$ and its fourth derivative has a discontinuity at $x=\ell/2$. Therefore, we have just $\mathcal{U}_0\in\mathcal{D}(A)$, but no greater regularity than that.

\begin{figure}
\begin{center}
    \begin{subfigure}[t]{0.5\textwidth}
        \centering
        \includegraphics[scale=0.50]{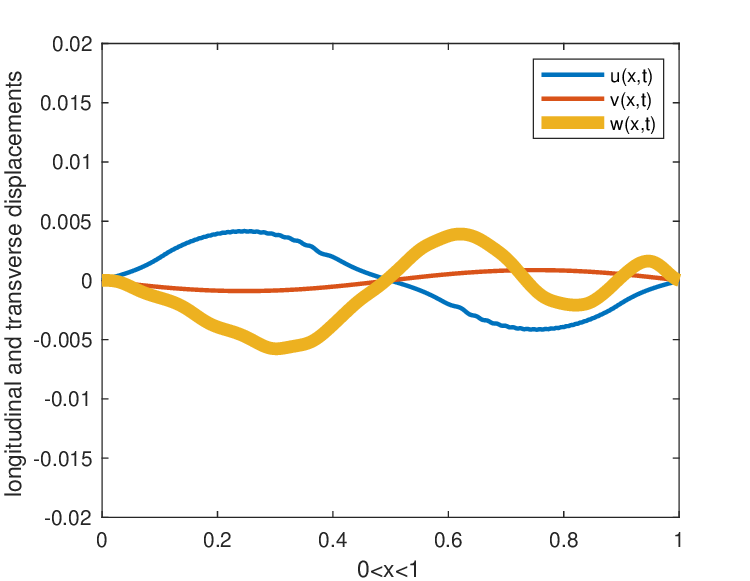}
        \caption{
State of the deformations at instant $t=10$}
    \end{subfigure}%
    \begin{subfigure}[t]{0.5\textwidth}
        \centering
        \includegraphics[scale=0.63]{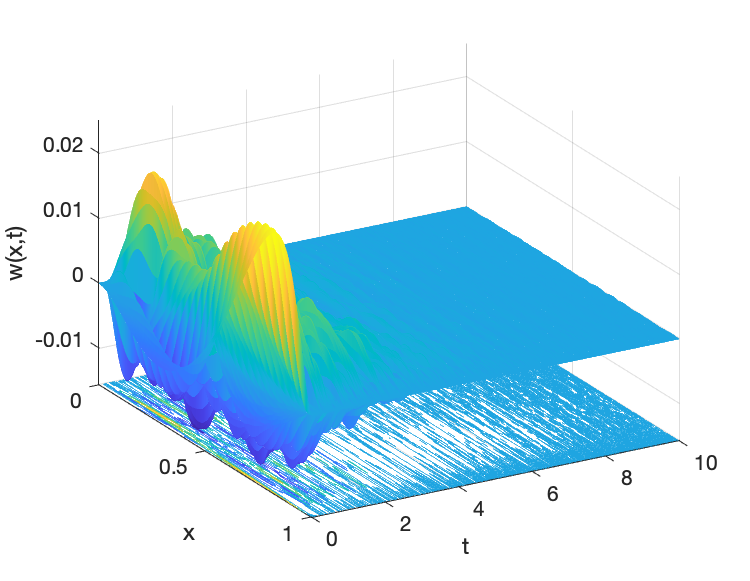} 
        \caption{Transversal deformations $w$}
    \end{subfigure}%
    \\
    \begin{subfigure}[t]{0.5\textwidth}
        \centering
        \includegraphics[scale=0.57]{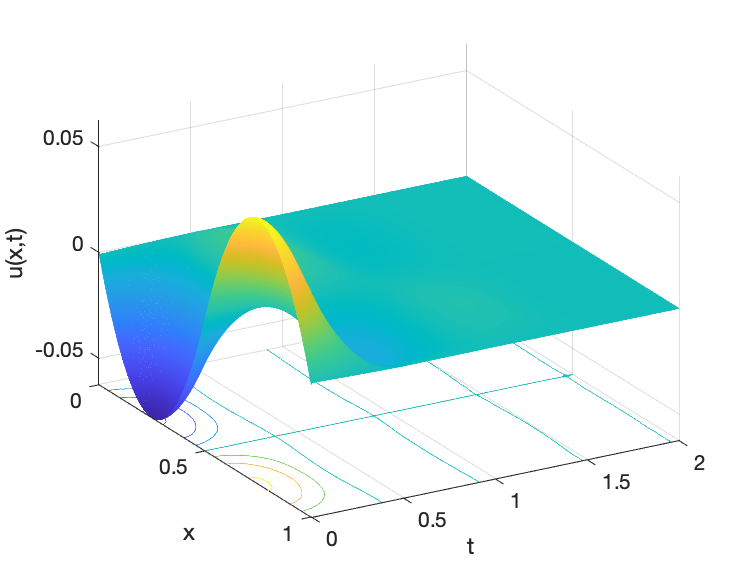} 
        \caption{Longitudinal deformations $u$}
    \end{subfigure}%
    \begin{subfigure}[t]{0.5\textwidth}
        \centering
        \includegraphics[scale=0.57]{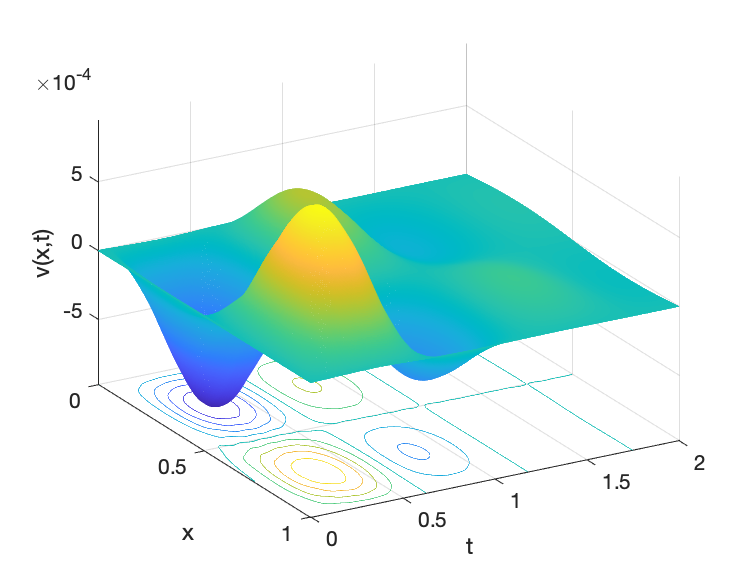} 
        \caption{Longitudinal deformations $v$}
    \end{subfigure}%
\caption{Simulation of the displacements. On the left: the displacements at $t=2$. On the right:
the evolution in time of the transverse displacements.
}\label{Fig2}
\end{center}
\end{figure}

In Figure \ref{Fig2}, we consider a fractional derivative  with $\alpha=0.5$,
and we do a simulation of the numerical scheme \eqref{702}, \eqref{703}, \eqref{Unmas1} with $J=500$, and 
$\delta=T/N_t$ with $T=100$ and $N_t=1000$.
The graphs in (A) of Figure \ref{Fig2} shows the state of
the instant of the deformations for $t=10$. The beam is represented by the thick yellow graph and its deformation corresponds to the transversal displacements. The longitudinal displacements of $u(x,10)$ and $v(x,10)$, are represented by the blue and red lines respectively. The graph (B) of the figure
\ref{Fig2} shows the evolution of the transverse displacements 
$w$ in term of space and time. Graphs (C) and (D) of the same figure show the longitudinal displacements $u$ and $v$ (respectively), also as a function of space and time.
In this set of figures, the decay of the quantities, and therefore of the energy, is visually observed due to the dispative terms of the fractional derivative type.

\begin{figure}
\begin{center}
\includegraphics[scale=0.57]{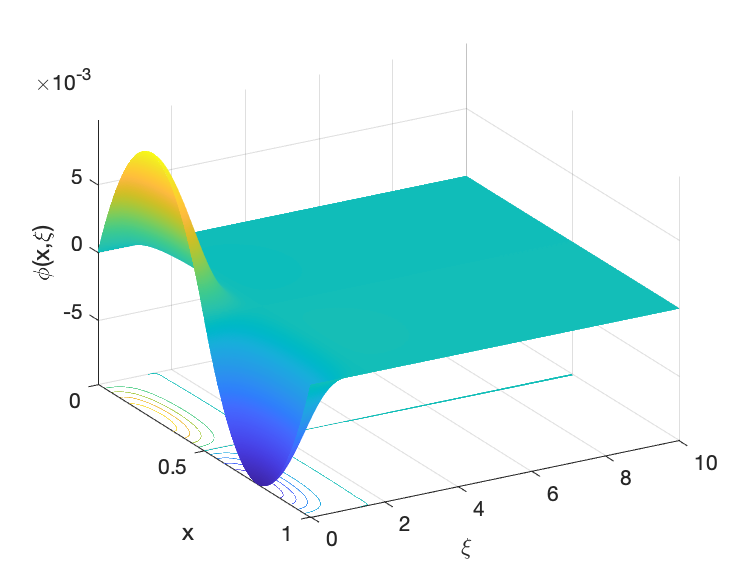}
\includegraphics[scale=0.57]
{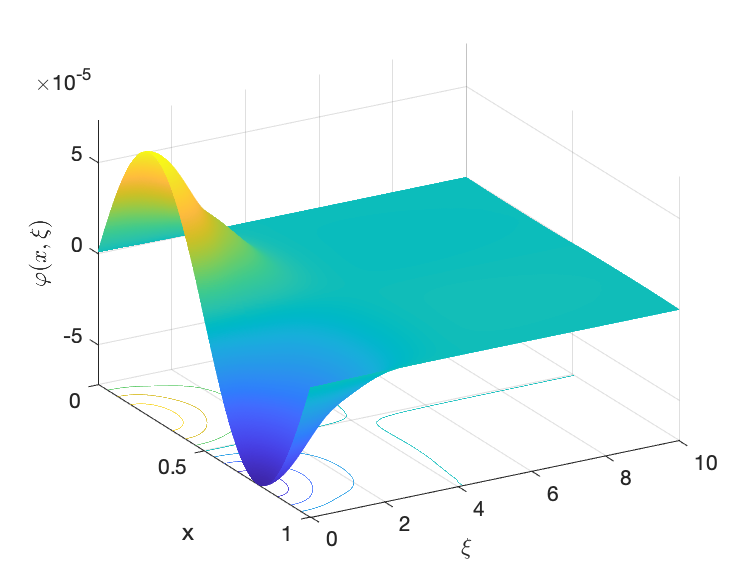} 
\includegraphics[scale=0.55]{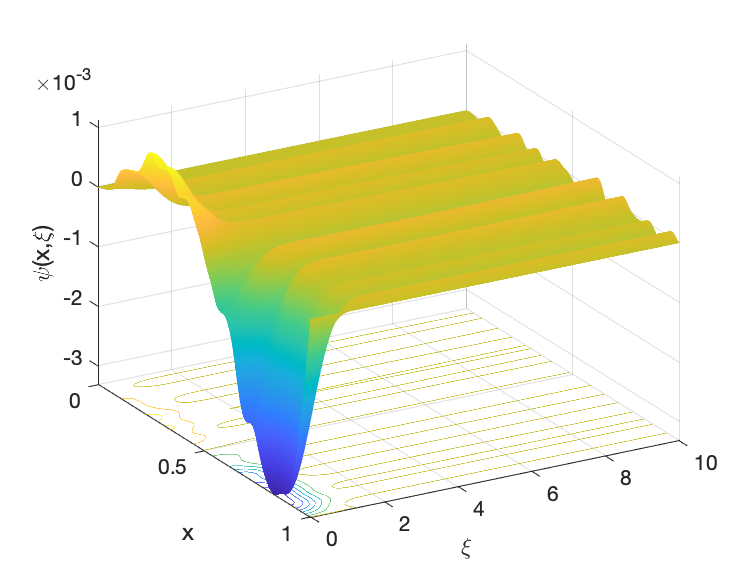} 
\caption{Behavior of the auxiliary variables $\phi$, $\varphi$ and $\psi$ of the augmented model \eqref{213}
}\label{Fig3}
\end{center}
\end{figure}
\medskip
On the other hand, we show in figure \ref{Fig3} the graphs of the numerical variables that approximate the auxiliary functions $\phi(x,t;\xi)$, $\varphi(x,t;\xi)$, and $\psi(x,t;\xi)$ of the augmented model \eqref{213}, for the instant $t=2$ of the simulations with the same parameters and data as the example in figure~\ref{Fig2}. We numerically observe a decay of these functions with respect to the variable $\xi$, although not so evident in the case of the variable $\psi$ associated with the term of the fractional derivative of the transverse displacement (when truncating at $\xi=10$). In this sense, we must be careful with the numerical truncation of the variable $\xi\in (0,+\infty)$. For this reason we will consider a truncation at $\xi=10,000$ for the examples that will follow, and that will show the polynomial decay of energy for different values of $\alpha$ and $\eta$.

\medskip


\begin{figure}
\begin{center}
    \begin{subfigure}[t]{0.5\textwidth}
        \centering
        \includegraphics[scale=0.61]{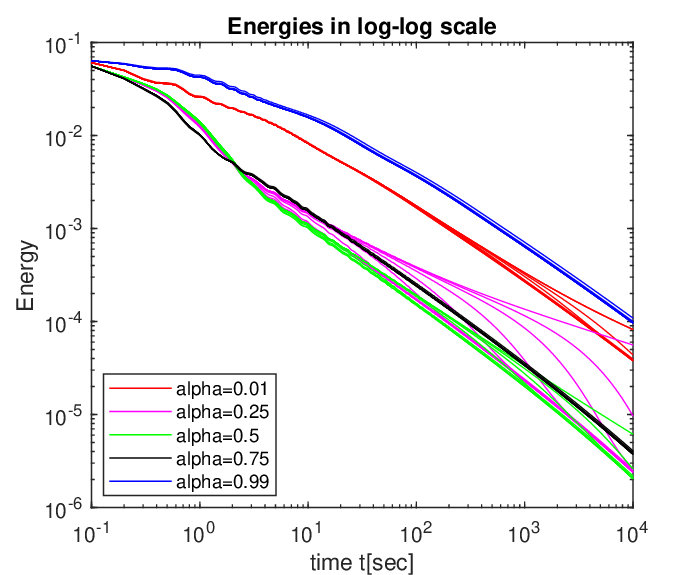} 
        \caption{Energies for different parameters  $\alpha$ and $\eta$\\ regrouped in terms of the values of $\alpha$.}
    \end{subfigure}%
    \begin{subfigure}[t]{0.5\textwidth}
        \centering
        \includegraphics[scale=0.55]{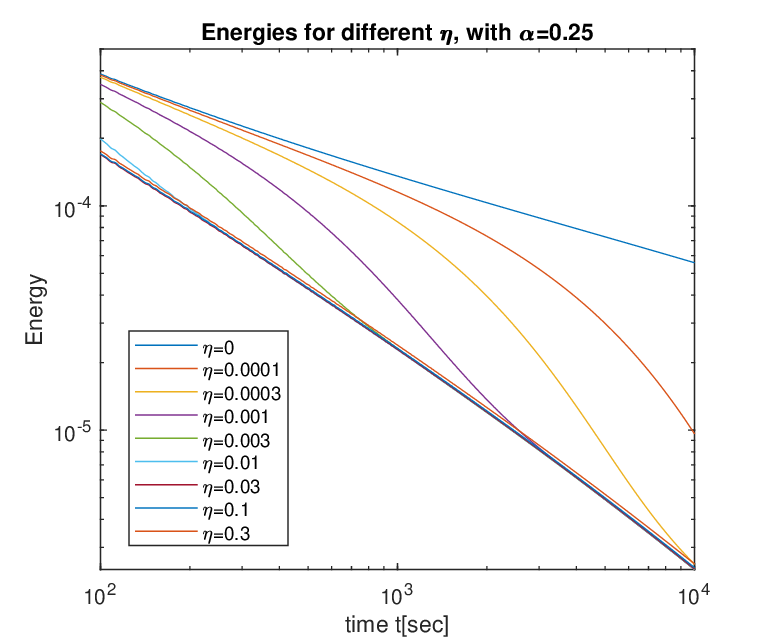} 
        \caption{Zoom: Energies for $\alpha=0.25$, \\ and different values of $\eta$}
    \end{subfigure}%
\\
    \begin{subfigure}[t]{0.5\textwidth}
        \centering
        \includegraphics[scale=0.55]{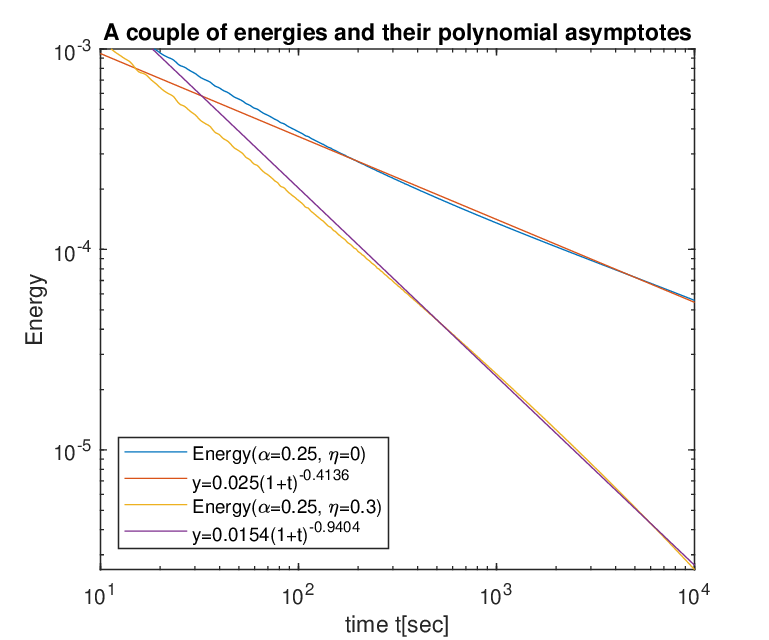} 
        \caption{Comparison between energy curves\\ with a same $\alpha=0.25$ as a function of time,\\ and their asymptotes $C(t+1)^p$ }
    \end{subfigure}%
    \begin{subfigure}[t]{0.5\textwidth}
        \centering
        \includegraphics[scale=0.55]{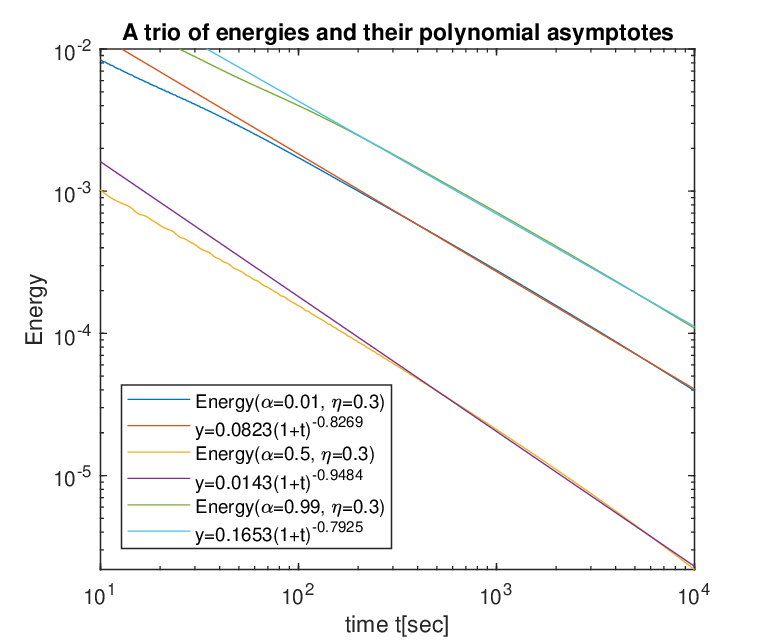} 
        \caption{Comparison between energy curves\\ with a same 
        $\eta=0.3$ as a function of time,\\ and their asymptotes $C(t+1)^p$ }
    \end{subfigure}%
\caption{
Energy with  different values of $\eta$ 
and $\alpha$ in loglog
scale.
}
\end{center}
\label{Fig4}
\end{figure}

\begin{table}
\begin{center}
    \begin{tabular}{cc|ccccc}
   &        & \multicolumn{5}{c}{\Large $\alpha$} \\
  \multicolumn{2}{c|}{$C(t+1)^p$} 
    & 0.01 & 0.25 & 0.50 & 0.75 & 0.99 \\
        \hline 
  \multirow{9}{*}{\Large $\eta$} &  $0$ & 
  \def\arraystretch{0.7}
  $\begin{array}{l} C=\ 0.0300 \\ p=\ -0.6465 \end{array}$ & 
  \def\arraystretch{0.7}
  $\begin{array}{l} 0.0250 \\ -0.4136 \end{array}$ & 
  \def\arraystretch{0.7}
  $\begin{array}{l} 0.0600 \\ -0.7497 \end{array}$ & 
  \def\arraystretch{0.7}
  $\begin{array}{l} 0.0186 \\ -0.9147 \end{array}$ & 
  \def\arraystretch{0.7}
  $\begin{array}{l} 0.1521 \\ -0.7968 \end{array}$
\\
\cline{2-7}
  &  $1.10^{-4}$ &  
  \def\arraystretch{0.7}
  $\begin{array}{l} C=\ 0.0869 \\ p=\ -0.8136 \end{array}$ & 
  \def\arraystretch{0.7}
  $\begin{array}{l} 0.0408 \\ -0.8630 \end{array}$ & 
   \def\arraystretch{0.7}
  $\begin{array}{l} 0.0250 \\ -0.9773 \end{array}$ &   
   \def\arraystretch{0.7}
  $\begin{array}{l} 0.0198 \\ -0.9255 \end{array}$ &
     \def\arraystretch{0.7}
  $\begin{array}{l} 0.1521 \\ -0.7968 \end{array}$   
  \\\cline{2-7}
  &  $3.10^{-4}$ &  
   \def\arraystretch{0.7}
   $\begin{array}{l} C=\ 0.1089 \\ p=\ -0.8600 \end{array}$ & 
  \def\arraystretch{0.7}
  $\begin{array}{l} 0.4175 \\  -1.2750 \end{array}$ & 
   \def\arraystretch{0.7}
  $\begin{array}{l} 0.0338 \\ -1.0469 \end{array}$ &   
   \def\arraystretch{0.7}
  $\begin{array}{l} 0.0198 \\ 0.9266 \end{array}$ &
     \def\arraystretch{0.7}
  $\begin{array}{l} 0.1521 \\ -0.7968 \end{array}$   
  \\\cline{2-7} 
   & $1.10^{-3}$ &  
      \def\arraystretch{0.7}
   $\begin{array}{l} C=\ 0.0894 \\ p=\ -0.8413 \end{array}$ & 
  \def\arraystretch{0.7}
  $\begin{array}{l} 0.0832 \\ -1.1394 \end{array}$ & 
   \def\arraystretch{0.7}
  $\begin{array}{l} 0.0338 \\ -0.9985 \end{array}$ &   
   \def\arraystretch{0.7}
  $\begin{array}{l} 0.0198 \\ -0.9236 \end{array}$ &
     \def\arraystretch{0.7}
  $\begin{array}{l} 0.1521 \\ -0.7968 \end{array}$   
  \\\cline{2-7}  
     & $3.10^{-3}$ &  
      \def\arraystretch{0.7}
   $\begin{array}{l} C=\ 0.0822 \\ p=\ -0.8317 \end{array}$ & 
  \def\arraystretch{0.7}
  $\begin{array}{l} 0.0258 \\ -1.0075 \end{array}$ & 
   \def\arraystretch{0.7}
  $\begin{array}{l} 0.0159 \\ -0.9688 \end{array}$ &   
   \def\arraystretch{0.7}
  $\begin{array}{l} 0.0190 \\ -0.9217 \end{array}$ &
     \def\arraystretch{0.7}
  $\begin{array}{l} 0.1522 \\ -0.7968 \end{array}$   
  \\\cline{2-7}  
      & $1.10^{-2}$ &  
      \def\arraystretch{0.7}
   $\begin{array}{l} C=\ 0.0803 \\ p=\ -0.8290 \end{array}$ & 
  \def\arraystretch{0.7}
  $\begin{array}{l} 0.0170 \\ -0.9586 \end{array}$ & 
   \def\arraystretch{0.7}
  $\begin{array}{l} 0.0144 \\ -0.9569 \end{array}$ &   
   \def\arraystretch{0.7}
  $\begin{array}{l} 0.0188 \\ -0.9207 \end{array}$ &
     \def\arraystretch{0.7}
  $\begin{array}{l} 0.1525 \\ -0.7967 \end{array}$   
  \\\cline{2-7}  
        & $3.10^{-2}$ &  
      \def\arraystretch{0.7}
   $\begin{array}{l} C=\ 0.0801 \\ p=\ -0.8283 \end{array}$ & 
  \def\arraystretch{0.7}
  $\begin{array}{l} 0.0155 \\ -0.9469 \end{array}$ & 
   \def\arraystretch{0.7}
  $\begin{array}{l} 0.0141 \\ -0.9534 \end{array}$ &   
   \def\arraystretch{0.7}
  $\begin{array}{l} 0.0188 \\ -0.9201 \end{array}$ &
     \def\arraystretch{0.7}
  $\begin{array}{l} 0.1533 \\ -0.7964 \end{array}$   
  \\\cline{2-7}  
          & $1.10^{-1}$ &  
      \def\arraystretch{0.7}
   $\begin{array}{l} C=\ 0.0806 \\ p=\ -0.8279 \end{array}$ & 
  \def\arraystretch{0.7}
  $\begin{array}{l} 0.0152 \\ -0.9431 \end{array}$ & 
   \def\arraystretch{0.7}
  $\begin{array}{l} 0.0140 \\ -0.9514 \end{array}$ &   
   \def\arraystretch{0.7}
  $\begin{array}{l} 0.0190 \\ -0.9190 \end{array}$ &
     \def\arraystretch{0.7}
  $\begin{array}{l} 0.1562 \\ -0.7955 \end{array}$   
   \\\cline{2-7}  
           & $3.10^{-1}$ &  
      \def\arraystretch{0.7}
   $\begin{array}{l} C=\ 0.0823 \\ p=\ -0.8269 \end{array}$ & 
  \def\arraystretch{0.7}
  $\begin{array}{l} 0.0154 \\ -0.9404 \end{array}$ & 
   \def\arraystretch{0.7}
  $\begin{array}{l} 0.0143 \\ -0.9485 \end{array}$ &   
   \def\arraystretch{0.7}
  $\begin{array}{l} 0.0196 \\ -0.9161 \end{array}$ &
     \def\arraystretch{0.7}
  $\begin{array}{l} 0.1653 \\ -0.7925 \end{array}$    
    \end{tabular}
    \caption{Constant $C$ and rate $p$ of polynomial decay $C(t+1)^p$ for the energy  by a least square method, for different values of 
    the fractional derivative parameters $\alpha$ and $\eta$.}
\label{Tab1}
\end{center}
\end{table}

\noindent
Finally, in order to numerically appreciate the polynomial decay of the energy, we do the simulation with the same parameters $\vartheta=\chi=\zeta=k=\gamma=1$, a beam of length $\ell=1$,
and the initial conditions \eqref{CI}, but for a longer time $T=10000$, and a time step $\delta t=0.1$. Likewise, in order to obtain better precision in the results, we truncate the variable $\xi$ of the auxiliary functions at $\xi=10000$.

In Figure 4, various plots depict energy on a log-log scale, where curves of the form $C(t+1)^p$ exhibit nearly linear behavior over  sufficiently long times.

The graph in Figure 4(A) illustrates energy curves for 9 different values of $\eta$ and 5 values of $\alpha$, totaling 45 combinations. Similar curves with the same $\alpha$ are grouped into 5 distinct colors. While there are variations in slopes $p$ and constants $C$, graphical observation reveals a close relationship between these parameters, particularly in their bounded ranges. The numerical estimation of $C$ and $p$ using a least squares method is presented in Table 1.

Additionally, from the graph, it is conjectured that energy decays more rapidly to zero for $\alpha=0.5$ (or similar values) and certain $\eta \neq 0$. Conversely, $\alpha=0.25$ exhibits intriguing behavior, accentuating differences across various $\eta$ values. Thus, the curves corresponding to $\alpha=0.25$ are separately plotted in Figure 4(B). Here, a pattern emerges with curves lying between two quasi-linear lines with distinct slopes. The steepest negative slope corresponds to $\eta=0.3$, while the smallest absolute slope corresponds to $\eta=0$, with intermediate values approaching the former asymptotically.

Figures 4(C) and 4(D) compare energy curves with their respective asymptotic forms resembling straight lines in a log-log graph. These asymptotic curves, of the form $C(1+t)^p$, are derived through least squares, effectively capturing the data's behavior. In Figure 4(C), distinct $\eta$ values for the same $\alpha$ exhibit notably different slopes (i.e., decay rates $p$). Conversely, in Figure 4(D), variations in the multiplicative constant $C$ outweigh changes in the decay rate $p$ across different $\alpha$ values.

\section*{Discussion and Conclusions}

We have demonstrated the existence, uniqueness, and stability of a Rao--Nakra model with dissipative terms of the generalized Caputo fractional derivative type within the domain, building upon ideas from a proven case for a broader class of evolution systems \cite{kais}. While focusing on the Rao--Nakra beam model, our treatment of the fractional term through an input--output relationship with a diffusion equation \cite{15} has illuminated a numerical method of finite-conservative differences type that ensures stable energy decay, unlike standard methods for approximating fractional derivatives. This enabled us to numerically simulate energy decay, confirming its polynomial nature.

One key observation from the data and interpretation of Figure 4 and Table 1 is the clear asymptotic polynomial decay behavior of energy. We note that the decay rate is notably greater in absolute value when $\eta \neq 0$, and lower when $\eta = 0$, validating the theoretical propositions (Proposition 4.3 and Theorem 4.2). However, we acknowledge discrepancies between the theoretically derived rates and the numerically estimated ones. We attribute this to factors such as the discretization step sizes $\delta t$ and $\delta x$, the order of convergence of the method introducing additional numerical dissipation, and the limited simulation times for observing true theoretical decay rates.

Future work should delve into numerical methods enhancing convergence order to obtain numerical values closer to theoretical ones. Additionally, exploring inverse or control problems, such as identifying parameter pairs $(\alpha, \eta)$ for quicker energy decay, merits attention. Notably, from Table 1, we observe a minimum in the rate $p$ near $(\alpha, \eta) = (0.5; 0.0003)$. An inverse problem could involve minimizing this slope with respect to $(\alpha, \eta)$. Moreover, considering the importance of the constant $C$, as evidenced by the curve nearest to zero at $t = 10,000$ corresponding to $(\alpha, \eta) = (0.5; 0.3)$ in Figure 4(A), remains intriguing.

However, beyond $t = 10,000 [s]$, maintaining a constant $p$ to ensure optimal decay is unclear. While minimizing energy decay speed is of interest, defining a suitable cost function characterizing this decay speed requires further exploration. Thus, identifying the appropriate cost function to minimize for solving the inverse problem remains an open task.

\bibliographystyle{elsarticle-num}

\end{document}